\documentclass[12pt]{article}
\usepackage{amsmath,amsfonts,amssymb,amsthm}

\usepackage[usenames,dvipsnames]{pstricks}
\usepackage{epsfig}
\usepackage{pst-grad} 
\usepackage{pst-plot}

\newcommand{\TG}{\mathcal{TG}}
\newcommand{\dih}[2]{DIH_{#1}(#2)}

\def\noi{\noindent}

\def\pf{\noi{\bf Proof.\ \,}}
\def\eop{{$\square$}}
\def\qed{{$\square$}}

\def\labtt#1{\label {#1} }
\def\labttr#1{\label {#1} }

\def\refpp#1{(\ref {#1})}

\def\a{\alpha}
\def\b{\beta}
\def\g{\gamma}
\def\d{\delta}

\def\s{\sigma}
\def\t{\tau}

\def\FF{{\mathbb F}}
\def\QQ{{\mathbb Q}}

\def\ZZ{{\mathbb Z}}

\def\la{\langle}
\def\ra{\rangle}
\def\<{\langle}
\def\>{\rangle}

\def\bs{\it}            

\def\dim{{\bs dim}}
\def\ker{{\bs Ker}}

\def\det{{\bs det}}

\def\half{{1 \over 2}}

\def\dual#1{#1^*}        

\def\dg#1{{\cal D}({#1})}

\begin{document}

\newtheorem{thm}{Theorem}[section]
\newtheorem{prop}[thm]{Proposition}
\newtheorem{lem}[thm]{Lemma}
\newtheorem{rem}[thm]{Remark}
\newtheorem{coro}[thm]{Corollary}
\newtheorem{conj}[thm]{Conjecture}
\newtheorem{de}[thm]{Definition}
\newtheorem{hyp}[thm]{Hypothesis}

\newtheorem{nota}[thm]{Notation}
\newtheorem{ex}[thm]{Example}
\newtheorem{proc}[thm]{Procedure}

\begin{center}\end{center}

\centerline{ \today  }
\begin{center}
{\Large
Diagonal lattices and rootless $EE_8$ pairs   }

\bigskip

\vspace{10mm}
Robert L.~Griess Jr.
\\[0pt]
Department of Mathematics\\[0pt] University of Michigan\\[0pt]
Ann Arbor, MI 48109  USA  \\[0pt]
{\tt rlg@umich.edu}\\[0pt]
\vskip 0.6cm

Ching Hung Lam
\\[0pt]
Institute of Mathematics \\[0pt]
Academia Sinica\\[0pt]
Taipei 10617, Taiwan\\[0pt]
{\tt chlam@math.sinica.edu.tw}\\[0pt]
\vskip 0.6cm
\end{center}

\begin{abstract}
Let $E$ be an integral lattice. We first discuss some general properties of an
SDC lattice, i.e., a sum of two diagonal copies of $E$  in $E\perp E$.  In
particular, we show that its group of isometries contains a wreath product. We
then specialize this study to the case of $E=E_8$ and provide a new and fairly
natural model for those rootless lattices which are sums of a pair of
$EE_8$-lattices.  This family of lattices was classified in \cite{GLEE8}.  We
prove that this set of isometry types is in bijection with the set of conjugacy
classes of rootless elements in the isometry group $O(E_8)$, i.e., those $h \in O(E_8)$ such that
the sublattice $(h-1)E_8$ contains no roots. Finally, our model gives new
embeddings of several of these lattices in the Leech lattice.
\end{abstract}

\bigskip 


\smallskip 


     \newpage
\tableofcontents

\section{Introduction}

In this article, {\it lattice} means a finitely generated free abelian group with a
rational valued
symmetric bilinear form.

We begin by defining the main construction used in this article.

\begin{nota}\labttr{sdcsetup}
Suppose that we are given an integral lattice, $E$, and an isometry $h\in O(E)$.
In $E\perp E$, we have two sublattices
\[
M:=\{(x,x) \mid x\in E\} \quad \text{ and }\quad  N:=\{(x, hx)\mid x \in E\}.
\]
Clearly, $M\cong N\cong \sqrt 2 E$ (where $\cong$ indicates isometry of
quadratic spaces). Define $L:=L(E,h):=M+N$.   We call $L$ an {\it SDC-lattice}
or, more precisely, an  $SCD(E,h)$-{\it lattice}  or $SDC(E)${\it -lattice},
meaning a sum of diagonal copies (of the fixed input lattice, $E$, using the
isometry $h$).
\end{nota}

Clearly, $L$ is integral (since it is a sublattice of $E\perp E$) and even (since the
generating set $M\cup N$ has only even norm vectors). Our first main result
shows that  $L$ has a large group of isometries \refpp{mainth1}.

\begin{thm}\labtt{mainth1}
Let $L, h$ be as in \refpp{sdcsetup}, where $h$ has order $n$.
Then $O(L)$ contains a chain of subgroups $\la
t_M, t_N\ra \le W_{M,N} \cong \ZZ_n \wr \ZZ_2$.
Furthermore, each of $t_M, t_N$ is a wreathing involution of $W_{M,N}$.
\end{thm}

Some lattices of great interest have this form. One has for instance the
Barnes-Wall lattices (for which $M, N$ are scaled copies of smaller rank
Barnes-Wall lattices and $h^2=-1$).    Additional examples are listed in
Section \ref{sec:2.1}.  One should note the trivial cases $h=1$, for which
$M=N$, and $h=-1$, for which $M+N=M\perp N$.

The term  {\it  $EE_8$-lattice} means a lattice isometric to $\sqrt 2 E_8$ \cite{GLEE8}.

We now consider rootless integral lattices spanned by a pair of $EE_8$-lattices.
They were studied and classified in \cite{GLEE8}.  Recently, we realized that they may be expressed as SDC-lattices \refpp{mainth2}.
The next two main results shows how
they may be expressed as SDC-lattices \refpp{mainth2}.

\begin{thm}\labtt{mainth2}
All rootless $EE_8$ pairs listed in \cite[Table 1]{GLEE8} can be embedded into
$E_8\perp E_8$ as $SDC(E_8)$-lattices  \refpp{sdcsetup}.
\end{thm}

\begin{thm}\labtt{mainth3}
There is a bijection between the conjugacy classes of rootless
elements in $O(E_8)$ and the isometry classes of rootless $EE_8$ pairs.
\end{thm}

An application  of modeling the lattices of \cite{GLEE8} as $SDC(E_8)$-lattices
is that one can see relatively natural embeddings of some of them into the
Leech lattice; see Section \ref{embed}.
 Such embeddings were first
demonstrated in \cite{GLEE8}, but the proofs were rather technical.

\bigskip

{\bf Conventions.  }
Group actions will be on the left.
Notations are generally standard.  We mention the relatively new notations $EE_8$ for $\sqrt 2 E_8$ \cite{GLEE8}, RSSD and SSD \refpp{rssd}.  For background on groups and lattices, see \cite{gal}.

\section{About SDC lattices}

In this section, $E$ is an arbitrary integral lattice.    Later in this article, we shall
specialize to the case $E=E_8$.

\medskip

\begin{de}\labtt{rssd}  A sublattice $X$ of an integral lattice $Y$ is called {\it RSSD}
if $2Y\le X+ann(X)$.  If $X$ is RSSD, the orthogonal transformation $t_X$ which
is $-1$ on $X$ and 1 on $ann(X)$ takes $Y$ to itself, whence $t_X\in O(Y)$.

The lattice $X$ is called {\it SSD} if $2\dual X\le X$.
An SSD lattice $X$ contained in the integral lattice $Y$ is RSSD in $Y$.
See \cite{GrE8,GLEE8,gal}.
\end{de}
\medskip

We use the notations of \refpp{sdcsetup}.

\begin{lem}\labtt{sdc1}   As maps on $E\perp E$, 
$t_M: (x,y)\mapsto (-y,-x)$ and
$t_N: (x,y)\mapsto (-h^{-1}y, -hx)$.
\end{lem}
\pf Direct calculation.  
Here is an argument for $t_M$. Write  $(x,y)=(\half
(x+y),\half (x+y)) + (\half (x-y), -\half (x-y))$ and note that the first summand on the right side is in $M$ and the second is in $ann(M)$.  Therefore, $t_M$ negates the first summand and fixes the second.

To verify the formula for $t_N$, notice that this map negates $N$ and fixes all
$(w,-hw)$, $w \in L$.  Then use the decomposition

$(x,y)=(\half (x+h^{-1}y), \half (hx+y)) + (\half (x-h^{-1}y), \half (-hx+y))$.
\eop

\begin{nota}\labtt{sdc2}
Define sublattices  $N':=\{ (x,h^{-1}x) \mid x \in E \}$ and $L':=M+N'$.

Define the following elements of $O(E\perp E)$:

$\b : (x,y)\mapsto (hx,y)$;

$\g : (x,y)\mapsto (x,hy)$;

$\d : (x,y) \mapsto (hx, hy)$;

$\d' : (x,y) \mapsto (h^{-1}x, hy)$.

\smallskip
\noindent
These maps satisfy
$\d = \b \g = \g \b$ and $\d' = \b^{-1} \g = \g \b^{-1}$.

We denote by $W(E,h)$ the group $\la t_M,  t_N, \b ,\g \ra$.  It is a subgroup of $O(E\perp E)$ (but we shall see that it embeds in $O(L)$ \refpp{sdc10}).  
\end{nota}

\begin{lem}\labtt{sdc3}  (i) $t_Nt_M=\d'$;

(ii) $\b = t_M\g t_M =  t_N\g t_N$;

(iii)  $W(E,h)$ is generated by any three of $t_M,  t_N, \b ,\g$.  Furthermore,
$W(E,h) = (\la \b \ra \times \la \g \ra)\la t_M \ra$ is isomorphic to the wreath product $\ZZ_{|h|} \wr \ZZ_2$;

(iv) $\la \b, \g \ra$ contains $\la \d, \d' \ra$ with index $(2,|h|)$.

(v) In $W(E,h)$, the stabilizer of $M$ is $\la t_M \ra \times \la \d \ra$ and the stabilizer of $N$ is $\la t_N \ra \times \la \d \ra$.

\end{lem}
\pf
(i) Direct calculation.

(ii) One may check the first equality by direct calculation.  For the second, note that $t_N=\d' t_M = t_M (\d')^{-1}$ and that $\d'$ and $\g$ commute.

(iii) Let $V$ be the subgroup of $W(E,h)$ generated by three of the generators
and let $H:=\la \b \ra \times \la \g \ra$. Then $V$ covers $W(E,h)/H\cong 2$,
i.e., $W(E,h)=HV$. If $V$ includes generators $\b, \g$, then $V\ge H$ and we
are done. If not, $V$ contains both $t_M$ and $t_N$, whence also $\d'$.
Clearly, $H$ is generated by any two of $\b, \g, \d'$ and so we conclude that
$V=W(E,h)$.

(iv) Clearly, $\la \b, \g \ra$ contains $\la \d, \d' \ra$.  The latter equals
$\la \b^2, \g^2, \d \ra$ and has index $(2,|h|)$ in $\la \b, \g \ra$.

(v) Let $S$ be the stabilizer of $M$ in $W(E,h)$.  We have  $\la t_M \ra \le S$.  
Since $W(E,h)= \la t_M \ra H$, the Dedekind law implies that 
$S=\la t_M \ra (S\cap H)$.  Clearly, $(S\cap H) = \la \d \ra$.  This completes the analysis for $M$.  The argument  for $N$ is similar.  
\eop

\begin{lem}\labtt{sdc4}   $\g (M)=N$,  $\g (N')=M$ and $\g (L')=L$.
\end{lem}

\begin{lem}\labtt{sdc5}  (i) $2L\le M+ann(M)$;

(ii) $2L' \le M+ann(M)$.
\end{lem}
\pf (i) It suffices to prove that $2N \le M+ann(M)$. An element of $N$ has shape
$(x,hx)$ for some $x\in E$.  We have $2(x,hx)=(x+hx,x+hx)+(x-hx,-x+hx)$. The
first summand is in $M$ and the second is in $ann(M)$.

(ii) Use (i) with $h$ replaced by $h^{-1}$. \eop

\begin{lem} \labtt{sdc5.5} $2L \le N+ann(N)$.
\end{lem}
\pf
Apply $\g$ to the containment \refpp{sdc5} (ii).

\begin{coro}\labtt{sdc6}  $\la t_M, t_N \ra$ maps $L$ to itself.  
\end{coro}
\pf  We have shown that $M$ and $N$ are RSSD lattices.  Therefore the
isometries $t_M$ and $t_N$ map $L$ to itself. \eop

\begin{rem}\labtt{sdc7}
The isometry group of $L$ contains an isomorphic copy $C(E,h)$ of $C_{O(E)}(h)$, acting
diagonally on $E \perp E$.
We have
 $\la -1, \d \ra \le C(E,h)$ and $C(E,h)$
 centralizes
$\la t_M, t_N \ra$.
\end{rem}

\begin{lem} \labtt{sdc9}
We have

(i) $L\cap (E\perp 0)= Im(h-1) \perp 0$; and

(ii) $L\cap (0\perp E)= 0 \perp Im(h-1)$.

\end{lem}
\pf (i)
Consider $a, b \in E$.  Then $(a,a)+(b,hb)\in E\perp 0$ if and only if $a=-hb$ if and only if $a+b=(1-h)b$.
This proves $L\cap (E\perp 0) \le Im(h-1) \perp 0$.
Conversely, suppose that
$c\in E$.  Then by \refpp{sdc1}, $((1-h)c,0)=(c,c)+(-hc,-c)=t_M(-(c,c)+(c,hc)) \in t_M(M+N)=M+N$ \refpp{sdc6}.   This proves $L\cap (E\perp 0) \ge Im(h-1) \perp 0$.

(ii) This follows from (i) and use of $t_M$ \refpp{sdc1}, \refpp{sdc6}. \eop

\begin{prop}\labtt{sdc10}
(i) 
$W(E,h)$ stabilizes $L$.  

(ii) The action of $W(E,h)$ on $L$ is faithful, so restriction gives an embedding of $W(E,h)$ in $O(L)$.  
\end{prop}
\pf (i) 
In view of \refpp{sdc3}(iii) and  \refpp{sdc6}, it suffices to prove that
$\g$ is in $O(L)$.
By \refpp{sdc4}, it suffices to prove that $\g (N)\le L$.
We take $a \in E$ and calculate
$\g (a,ha)=(a,h^2a)=(a,ha)+(0,h^2a-ha))$.
Obviously,  $(a, ha)\in N \le L$. We have
$(0,h^2a-ha))=(0,(h-1)ha)$, which is in $L\cap (0\perp E)$ by \refpp{sdc9}, so we are done.

(ii) Let $K$ be the kernel of the action of $W(E,h)$ on $L$.  We may 
assume that $E\ne 0$.  
By \refpp{sdc3}(v), $K \le \la t_M, \d \ra$.

We shall argue that $K \le \la \d \ra$.  Suppose otherwise.  Consider an integer $i$ so that $z:=\d^i t_M \in K$.  
Then $z$ takes $(x,x)$ to $(-h^ix, -h^ix)$ which is $(x,x)$ since $z\in K$.  It follows that $h^i=-1$ on $E$.  
By \refpp{sdc1}, $z$ takes $(x, hx)$ to $(hx,x)$, which must equal $(x,hx)$, for all $x \in E$.  
We   conclude that $h=1$.  
Since $E\ne 0$, this incompatible with $h^i = -1$.  

We have $K \le \la \d \ra$.  Since the group 
$\la \d \ra$ acts faithfully on $M$, it acts faithfully on $L$ and we conclude that $K=1$.   
\eop

\begin{lem}\labtt{annm(n)}
Let $M$ and $N$ be defined as above. Then
\[
\begin{split}
ann_N(M) =&\{ (\a,-\a)\mid \a\in E\text{ and } h\a =-\a\},\quad \text{  and  }\\
ann_M(N) =& \{ (\a,\a)\mid \a\in E\text{ and } h\a =-\a\} .
\end{split}
\]

\end{lem}

\pf
We prove the first equality.  The proof of the second is similar.

Let $(\a, h\a)\in N$.  Then
\[
\begin{split}
&( \a,h\a) \text{ annihilates } M  \\
\text{  if and only if }\quad   & (\a , \b) +( h\a,\b) =0 \text{ for all } \b\in E \\
\text{  if and only if }\quad  & (h\a  +\a, \b)  =0 \text{ for all } \b\in E \\
\text{  if and only if }\quad  & h\a=-\a.
\end{split}
\]
Thus, $ ann_N(M) =\{ (\a,-\a)\in E\perp E\mid \a\in E\text{ and } h\a=-\a\}$
as desired. \qed

\begin{rem}\labtt{sdc8}  (i) Given a pair of isometric
doubly  even lattices, $M, N$ in Euclidean space, such that $M+N$ is integral
and $M, N$ are RSSD in $M+N$, when is there a representation of $M+N$ in the
form of \refpp{sdcsetup}? One would need to define a suitable $h$.  The
following example indicates a caution.

Let the lattice $L$ have basis $u, v$ and Gram matrix
 $\begin{pmatrix}
 2a&b \cr
 b&2a
 \end{pmatrix}$,
for integers $a\ge 1$ and $b$.  For positive definiteness, we require 
$4a^2-b^2>0$.   The $A_2$-lattice is such an example.

Let $E$ be the rank 1 lattice with Gram matrix $(a)$.  Then $M:=span\{ u \}$ and
$N:=span\{ v \}$ are sublattices of $L$ isometric to $\sqrt 2 E$ and their sum is
$L$.  
The condition that $M$ and $N$ be RSSD in $L$ is $a|b$.  

If $L$ were isometric to $SDC(E,h)$ with $M, N$ as in \refpp{sdcsetup}, then $h=\pm 1$ and so 
$b\in 2a\ZZ$, which implies the RSSD condition $a|b$.  
The necessary condition $b\in 2a\ZZ$ implies that $L$ is not positive definite if $b\ne 0$, so the above $L$ are not  $SDC(E,h)$ if $b\ne 0$.

(ii) A study of SDC lattices was carried out by Paul Lewis in his 2010
undergraduate research project \cite{lewis}. For many cases of familiar input lattice
$E$ and isometry $h$, the resulting $SDC(E, h)$ is another familiar lattice,
but there are surprises.
\end{rem}

\section{About rootless isometries}

We continue to use the notations \refpp{sdcsetup}.

\begin{de}
We say $h\in O(E)$ is rootless if $(h-1)E$ contains no roots.
\end{de}

\begin{lem}\labtt{M+N(h-1)L}  Let $E$ be an even lattice.
The sum $M+N$ is rootless if and only if $h$ is rootless.
\end{lem}

\pf  Let $x= (\a+\b, \a+h\b)\in M+N$, where $\a,\b\in E$. If both $\a+h\b$ and
$\a+\b$ are non-zero, then  $(x,x) \geq 2+2=4$.

If $\a+\b=0$, then $x=(0,(h-1)\b)$ and if $\a+h\b=0$, then $x=(-(h-1)\b, 0)$.
Thus, $(x,x) > 2$ if $(h-1)E$ is rootless.

On the other hand,  $(0,(h-1)\a)\in M+N $ for any $\a\in E$.  Therefore,
$(h-1)E$ is rootless if $M+N$ is. \qed

\medskip

We now take $E$ to be $E_8$ and begin determination of those $h$ for which
the conditions of \refpp{M+N(h-1)L} hold.

\begin{lem}
Suppose that $h\in O(E)$ and $h$ is rootless.  Then so is $h^i$ for all $i\in \ZZ$.
\end{lem}
\pf  We may assume that $i \ge 1$.
Since $h^i-1=(h-1)(1+h+h^2+\dots +h^{i-1})$,
this is clear.
\eop

\smallskip

\begin{nota}\labtt{ppart}
Recall that if $g$ is a group element of finite order $mn$, with $(m,n)=1$, then
$g$ is uniquely expressible as $g=hk$,  where $h$ has order $m$ and $k $ has
order $n$ and $hk=kh$.  Such $h, k$ lie in $\la g \ra$.   If $m$ is a power of the
prime $p$, we call $h$, $k$ the {\it $p$-part,  $p'$-part} of $g$, respectively.
Denote by $g_p, g_{p'}$ be the $p$-part, $p'$-part of $g$, respectively.
\end{nota}

\begin{coro}  If $h\in O(E)$ is rootless, then so are the $p$-parts of $h$, for all primes $p$.
\end{coro}

\begin{coro}
Suppose that $E$ contains roots, that $h\in O(E)$ is rootless and that $p, q$ are
distinct primes so that  $pq| |h|$.   Then at most one of  $h_p, h_q$ has no
eigenvalue 1.
\end{coro}
\pf If $h_p$ has no eigenvalue 1, $(h_p-1)E$ has index a power of $p$. If $h_q$
has no eigenvalue 1, $(h_q-1)E$ has index a power of $q$. If both of these
statements are true then $(h-1)E$ contains $(h_p-1)E + (h_q-1)E$, which by
relative primeness has index 1 in $E$.  This contradicts the rootless property of
$h$. \eop

\subsection{Root lattice of type $A$}
We shall review some basic properties of the root lattices of type $A_n$.

We use the {\it standard model} for $A_n$, i.e.,
\[
A_n=\left \{(x_1, x_2, \dots, x_{n+1}) \in \ZZ^{n+1} \, \left |\,  \sum_{i=1}^{n+1} x_i =0\right. \right \}.
\]
Then the roots of $A_n$ are given by
\[
\{\pm (e_i-e_j)\mid 1\leq i<j \leq {n+1}\},
\]
where $\{e_1=(1,0,\dots,0), \dots, e_{n+1} =(0, 0, \dots, 1)\}$ is the standard basis
of $\ZZ^{n+1}$.

\begin{nota}\labttr{gammaj}
Recall that $(A_n^*)/ A_n \cong \ZZ_{n+1}$. Let $\gamma_{A_n}(0)=0$ and
\[
\gamma_{A_n}(j) = \frac{1}{n+1} \left
(-(n+1-j)\sum_{i=1}^j e_i+ j \sum_{i=j+1}^{n+1} e_i\right),  \text{ for } j = 1, \dots ,
n.
\]
Then $\gamma_{A_n}(j) \in A_n^*$. In fact, $\{\gamma_{A_n}(0), \gamma_{A_n}(1),
\dots, \gamma_{A_n}(n)\}$ forms a transversal of $A_n$ in $A_n^*$ \cite[Chapter
4]{CS}. We also note that the norm of $\gamma_{A_n}(j)$ is equal to
$j(n+1-j)/(n+1)$ for all $j=0, \dots, n$.
\end{nota}

\begin{nota}
Let $h_{A_n}$ be an $(n+1)$-cycle  in $Weyl (A_n)\cong Sym_{n+1}$.
\end{nota}

\begin{lem} \label{h-1A*root}
For $j=1, \dots, n$,  $(h_{A_n}-1)(\gamma_{A_n}(j)) $ is a root.
\end{lem}

\pf By definition, $(h_{A_n}-1)(\gamma_{A_n}(j)) = e_1-e_{j+1}$ is a root. \eop

\begin{lem}\label{h-1An}
$( h_{A_n}-1) A_n$ is rootless.
\end{lem}

\pf We may assume $h_{A_n}$ is the cyclic permutation of the
$n+1$-coordinates.

Suppose $(h_{A_n}-1)\a  $ is a root for some $\a= (x_1, x_2, \dots, x_{n+1})\in
A_n$. Without loss, we may assume $ (h_{A_n}-1)\a= e_1-e_j$ for some $j\ge 2$.

Then we have
\[
x_{n+1}-x_{1} =1,\  x_{j-1}-x_{j}=-1,\   x_1=\dots=x_{j-1}  \text{  and  }   x_j=\dots=x_{n+1}.
\]
That implies $x_{n+1}=1+x_1$. Moreover, $x_1+\cdots+x_{n+1}=0$. Thus, we have
$(j-1)x_1+(n+2-j)(x_1+1)=0$ or $x_1=-\frac{n+2-j}{n+1}$, which is not an
integer since $2\le j \le n+1$, a contradiction. \qed

\begin{lem}\label{h-1An*}
Let $A_n^*$ be the dual lattice of $A_n$. Then $ (h_{A_n}-1)A_n^* = A_n$
\end{lem}

\pf
{\it First proof: }
Again, we shall use the standard model for $A_n$.  Then $A_n^* $ is the
$\ZZ$-span of
\[
\frac{1}{n+1} (1, 1,1,  \dots, 1, -n),  \frac{1}{n+1} (1, 1, \dots,-n,1), \dots,
\frac{1}{n+1}(1, -n, 1, \dots, 1,1).
\]

Note that
\[
(h_{A_n} -1)\left (\frac{1}{n+1} (1, 1,1,  \dots, 1, -n) \right)  = (1, 0,\dots, 0, -1) \in A_n.
\]
Similarly, we can show that $(h_{A_n}- 1)A_n^* \leq  A_n$.

On the other hand, the set $$\{(1,0,\dots, 0,-1),  (0,0,\dots, -1,1), \dots,  (0,-1,1,
\dots,0) \}$$ spans $A_n$ and hence  $ (h_{A_n}- 1)A_n^*= A_n$.

{\it Second proof: } Since $(h-1)\dual {A_n}=(h-1)\ZZ^{n+1}=span\{ e_i-e_{i+1}
\mid i=1,2 \dots \}$,  this is clear. \qed

\begin{lem}\labtt{aine8}  Let $X$ be a type $A_m$ lattice contained in $E_8$.  Then $X$ is a direct summand unless $m=8$.
\end{lem}
\pf If $X$ is properly contained in a summand, $S$, of $E_8$, then there exists
an integer $d\ge 2$ so that $d^2|det(X)$.  Since $\det(X)=m+1$ and $m \le 8$,
$m=3$ or $m=8$. If $m=3$, $d=2$ and so $det(S)=1$, whence $S\cong \ZZ^4$,
which is an odd lattice, a contradiction.   Therefore, $m=8$.
\eop

\begin{lem}\labtt{2non0}
Identify $Q:=A_{i_1} \perp \cdots\perp A_{i_\ell}$ with  a rank $8$ sublattices of $E_8$.
For any  $1\leq k\leq \ell$, define $h:=h_k:=h_{A_{i_1}}\oplus \cdots \oplus h_{
A_{i_k}} \oplus id\oplus  \cdots \oplus id$.

(a)  Suppose that for any $x\in E_8\setminus Q$,  $(h-1)x$  is either $0$ or has
non-zero projections to at least two of the $A_i$'s. Then $(h-1)E_8$ is rootless.

(b) Suppose there exists an element $x\in E_8\setminus Q$ such that $(h-1)x$
has non-zero projections to exactly one of the $A_i$'s. Then $(h-1)E_8$ has a
root.
\end{lem}

\pf  (a)  By Lemma \ref{h-1An},  it is clear that $(h-1)Q $ has no roots. Now let
$x\in E_8\setminus Q$. Then by our assumption and Lemma \ref{h-1An*},
$(h-1)x$  is either $0$ or has norm $\geq 2\times 2$. Hence, $(h-1)E_8$ has no
roots.

(b)   Let  $x\in E_8\setminus Q$ such that $(h-1)x$ has non-zero projections to
exactly one of the $A_i$'s, say to $A_{i_1}$.

Let $a$ be the projection of $x$ to $A_{i_1}^*$. Then there exists $j\in\{1,
\dots, i_1\}$ such that $a$ is in the coset $\gamma_{A_{i_1}}(j) +A_{i_1}$ (cf.
Notation \ref{gammaj}). Thus, there exists $b\in A_{i_1}$ such that $
a+b=\gamma_{A_{i_1}}(j)$. In this case,
\[
(h-1)(x+b) = (h_{A_{i_1}} -1)(a+b) =
(h_{A_{i_1}}-1) (\gamma_j),
\]
which is a root by Lemma \refpp{h-1A*root}. \qed

\medskip

\section{Eliminating cases}

We begin to study the cases where $h$ is $p$-element for some prime $p$. Recall
that $O(E_8)$ has order  $2^{14}{\cdot} 3^5{\cdot} 5^2{\cdot} 7$.

\smallskip

{\bf Convention.  } When we consider an embedding of lattices $X \le Y$, we may
describe it informally as  containment of isometry types, for example ``$A_1^8 \le
E_8$" or ``$A_2^3 \le E_6$".   Given such a containment, one may use notations
for isometries of the sublattice and make use of their unique extensions to
overlattices. This informally should not cause confusion.

\subsection{The prime 7}

\begin{lem}
There is no rootless element of order $7$ in $O(E_8)$.
\end{lem}
\pf By Sylow's theorem, there is only one conjugacy class of order $7$
subgroups in $O(E_8)$. Without loss, we may assume
\[
h=h_{A_6}\oplus id_B,
\]
where $B=ann_{E_8}(A_6)$. However, $(h-1)E_8$ has roots by Lemma
\ref{h-1A*root}. \eop

\subsection{The prime 5}

\begin{thm}\labtt{rootless5}
A rootless element of order $5$ is fixed point free and is
conjugate to $h_{A_4}\oplus h_{A_4}$.
\end{thm}

\pf  Let $h$ be an order 5 in $O(E_8)$. Then there is a root $\a$ such that
$h\a\neq \a$ since $E_8$ is generated by roots. Then, $(h^4+h^3+h^2+h+1)(\a)=0$
and $( (h^4+h^3+h^2+h+1)(\a), \a)=0$. This implies $( h\a, \a)+ ( h^2\a, \a)
=-1$ since $( h\a, \a)=( h^4\a, \a)$, $(h^2\a, \a)= (h^3\a, \a)$ and $(\a,
\a)=2$. By Cauchy-Schwarz inequality, we have $|( h\a, \a)|< 2$ and $|(h^2\a,
\a)|< 2$ and thus $( h\a, \a)=-1$ , $(h^2\a, \a)=0$ or $(h\a, \a)=0$, $( h^2\a,
\a) =-1$. Therefore, $K=span\{h^i\a\mid 0\leq i\leq 3\} \cong A_4$ since the
Gram matrix of $K$ is given by
\[
\begin{pmatrix}
2 &-1 &0&0\\
-1&2&-1&0\\
0&-1&2&-1\\
0&0&-1&2
\end{pmatrix}.
\]
Then $ann_{E_8}(K)\cong A_4$ \cite[(5.3.2)]{gal} and $h$ stabilizes both $K$
and $ann_{E_8}(K)$.

\textbf{Case 1:}  $h$ fixes $ann_{E_8}(K)$ pointwise. Then $h$ is conjugate to
$h_{A_4}\oplus id_{A_4}$, which is not rootless by \refpp{2non0} (b).

\textbf{Case 2:}  There exists a root $\b\in ann_{E_8}(K)$ such that $h\b\neq
\b$. Then $ann_{E_8}(K)= span\{h^i\b\mid 0\leq i\leq 4\}\cong A_4$. In this
case, $h$ is fixed point free and lies in $Weyl(K)\times Weyl( ann_{E_8}(K) )\cong
Sym_5 \times Sym_5$. Such elements form a single conjugacy class, so $h$ is
conjugate to $h_{A_4}\oplus h_{A_4}$ and $h$ is rootless \refpp{h-1An*}. \eop

\subsection{The prime 3}

\noindent \textbf{Order 3}

\begin{nota}\labtt{order3setup}
Let $h$ be an element of order 3 in $O(E_8)$. Let $F$ be the
fixed point sublattice of $h$ in $E_8$.
Let $J:=ann_{E_8}(F)$.
\end{nota}

By the analysis in \cite{GLEE8}, $\dg{F}\cong 3^s$ for some integer $s$. Thus,
by \cite[Lemma D.9]{GLEE8}, $F\cong 0, A_2, A_2\perp A_2,$ or $E_6$. Note that
in each case, $F$ contains an orthogonal direct sum of $A_2$'s with finite
index.

We have $J\cong E_8, E_6, A_2\perp A_2$ and $A_2$,
respectively and $h$ is fixed point free on $J$. Recall that the fixed point
free elements of order $3$ in $O(E_8)$ form one conjugacy class and they are
conjugate to $h_{A_2}^{\oplus 4}$ in $O(E_8)$. The fixed point free elements of
order $3$ also form one conjugacy class in $O(E_6)$ and they are conjugate to
$h_{A_2}^{\oplus 3}$ (see for example \cite{Atlas}). Therefore, in each case,
there exists a sublattice of $E_8$ which we may identify with  $A_2^4$ such
that $h=h_{A_2}^{\oplus 4-k}\oplus id_{A_2}^{\oplus k}$, where $k=\frac{1}2
\dim F$. Recall that $E_8/A_2^4$ can be identified with the tetracode
$\mathcal{C}_4$, which is a self-dual code of length $4$, minimal weight $3$
\cite{CS,gr12}.  Now, by Lemma \ref{2non0}, we have the theorem.

\begin{thm}\labtt{rootless3}
Let $h$ be an element of order 3 in $O(E_8)$. Then $h$ is rootless if and only
if $F=Fix(h)=0$ or $\cong A_2$. Identify $A_2^4$ with a sublattice of $E_8$.
Then,  $h$ is conjugate to $h_{A_2}^{\oplus 4}$ if $F=0$ and $h_{A_2}^{\oplus
3} \oplus id_{A_2}$ if $F\cong A_2$.
\end{thm}

\medskip

\noindent \textbf{Order 9}
\medskip

\begin{nota}\labtt{order9setup} Let $h$ be an element of order 9 in $O(E_8)$.
Let $g:=h^3$ and $F:=Fix(h^3)=\ker(h^3-1)$.  Let
$J:=ann_{E_8}(F)$.
\end{nota}

Then the minimal polynomial of $h$ on $J$ is divisible by the irreducible
cyclotomic polynomial $x^6+x^3+1$ and the minimal polynomial for $h$ on $F$
is $x-1$ or $x^2+x+1$. Hence,  $rank(F)=2$ (whence $F\cong A_2$) and $rank(
J)$ is 6.  Since $h$ stabilizes both $F$ and $J=ann_{E_8}(F)\cong E_6$, $h|_{F}$
defines an element of order 1 or 3 in $O(F)$ and $h|_{J}$ is an order 9 element in
$O(J)$.

\begin{lem}\labtt{zpwrzp}
In $\ZZ_p\wr \ZZ_p$, there are $(p-1)^2$ conjugacy classes of elements of order
$p^2$. More precisely, we let $B=B_1\times\dots \times B_p$ where each factor $B_i$ has
order $p$ and the order $p$ automorphism $g$ acts on $B$ by cyclically
permuting the $p$ factors. Thus the semidirect product $B\la g\ra$ is
isomorphic to $\ZZ_p\wr \ZZ_p$. The classes of order $p^2$ are represented by
${u_i}^kg^i$, $i=1,2,\dots , p-1$, $k=1, \dots , p-1$, where for each $i$,  $u_i$ is a generator for $B$ as
a $\la g\ra$-module.
\end{lem}
\pf  We count. Two such elements $u_i^kg^i$ and $u_j^{\ell} ig^j$ can not be
conjugate if $i \ne j$ or $k \ne \ell$ since their images modulo $(B\la g
\ra)'$ are distinct. The conjugacy class of such an element has cardinality
$p^{p-1}$ since $B$ is a free module for $B\la  g \ra/B$. Therefore, we have
accounted for $(p-1)(p^p - p^{p-1})$ elements of $B\la g \ra$. The $p^p$
elements of $B$ have order 1 or $p$. If $i=1,2,\dots , p-1$ and $v\in B$ does
not generate $B$ as a $\la g\ra$-module, then $vg^i$ has order 1 or $p$. This
latter category accounts for the remaining $(p-1)p^{p-1}$ elements of $B\la
g\ra$. \eop

\begin{coro}\labtt{order9ine6}
In $O(E_6)$, there is just one conjugacy class of elements of order 9.
\end{coro}
\pf
We view the $E_6$ lattice as an overlattice of $A_2^3$, defined by  glue
vector $(1,1,1)$.  From this viewpoint, it is obvious that we have a group of
automorphisms $H:=Weyl(A_2)\wr Sym_3$. The analysis of \refpp{zpwrzp} shows
that we have exactly four conjugacy classes of elements of order 9 in a Sylow
3-subgroup of $H$.  These classes are fused in a Sylow 3-normalizer in $H$.
\eop

\begin{thm}\labtt{order9}
There are no rootless elements of order 9 in $O(E_8)$.
\end{thm}

\pf
Let $h'= h|_{J}\in O(E_6)$ be an element of order 9. Recall that $E_6$
contains a sublattice of type $A_2^3$ and we may assume
\[
E_6 = span\{ A_2^3,  (\gamma, \gamma, \gamma)\}
\]
where $\gamma= \frac{1}3(1,1,-2)\in A_2^*$.

Note that there is only one conjugacy class of order 9 in $O(E_6)$
\refpp{order9ine6}. Thus, we may assume that $h'=\tau \sigma $, where $\sigma =
h_{A_2} \oplus id_{A_2} \oplus id_{A_2}$ and $\tau$ is  a cyclic permutation of
the 3 copies of $A_2$.

Let $ \a= (\gamma, \gamma, \gamma)= \frac{1}3( 1,1,-2; 1,1,-2; 1,1,-2).$ Then
\[
h'(\a)= \frac{1}3( 1,1,-2; -2,1,1; 1,1,-2) \text{ and } (h'-1) \a= (0,0,0;1,0,-1;0,0,0),
\]
which is a root. \eop

There are no elements of order 27 in $O(E_8)$, by \refpp{zpwrzp} and the fact
that $Weyl(E_6)\times Weyl(A_2)$ embeds with index prime to 3 in $O(E_8)$.
Therefore, we have treated all cases of 3-elements in $O(E_8)$.

\subsection{The prime 2}

{\bf Order 2}

\medskip

Suppose  $h\in O(E_8)$  has order 2. Then  the $(-1)$-eigenlattice  $L^-(h)$ of
$h$ is a RSSD sublattice of $E_8$. By the classification of RSSD lattices in $E_8$
\cite[Lemma D.2]{GLEE8}, there are nine possible cases up to conjugation and
$L^-(h) \cong  A_1^k, k\leq 4, D_4, D_4\perp A_1, D_6, E_7$ or $E_8$. For each
case, there exists a sublattice $A_1^8 <  E_8$ such that $h=h_{A_1}^{\oplus
k}\oplus id_{A_1}^{\oplus (8-k)}$, where $k=\dim L^-(h)$ (proof: each of the
above RSSD lattices contains an orthogonal direct sum of $A_1$s with finite
index).

\begin{thm}\label{order2}
 Suppose  $h\in O(E_8)$  has order 2. Then $h$ is rootless if and only if
 $L^-(h)\cong D_4, D_6, E_7$ or $E_8$.
\end{thm}

\pf Suppose $\dim (L^-(h))= k$. Then there exists $\a_1, \dots, \a_k\in L^-(h)$
such that $( \a_i,\a_j)=2\delta_{i,j}$ for $i,j=1, \dots,k$. Take $\a_{k+1},
\dots, \a_8\in ann_{E_8}(L^-(h))$ such that
\[
A=\ZZ\a_1\oplus\cdots \oplus \ZZ\a_8 \cong A_1^8.
\]
Then the quotient group $E_8/A$ can be identified  with the Hamming $[8,4,4]$
code $H_8$.

\textbf{Case 1:} $L^-(h)\cong A_1^k$, $1\leq k \leq 4$. By identifying  a
codeword with its support, we know that $\{1, \dots, k\}\notin H_8$ since the
minimal weight of $H_8$ is $4$ and $L^-(h)\cong D_4$ if $\{1,2,3,4\}\in H_8$.
Hence there exists $a\in H_8$ such that $|\{1, \dots, k\}\cap a|$ is odd. Without
loss, we may assume $a$ has weight $4$. Then $|\{1, \dots, k\}\cap a|=1$ or
$3$.

If $|\{1, \dots, k\}\cap a|=1$, let $\a_{a}= \frac{1}2 \sum_{i\in a}\a_i$. Then
$(h-1)\a_a = -\a_j$ is a root, where $\{ j \}= \{1, \dots, k\}\cap a$.
If $|\{1, \dots, k\}\cap a|=3$, let $\bar{a}= \{1, \dots, 8\}\setminus a$. Then
$|\{1, \dots, k\}\cap \bar{a}|=1$ and we get a contradiction as before.
We conclude that $h$ is not rootless.

\textbf{Case 2:} $L^-(h)\cong D_4\oplus A_1$. Then $k=5$. There exists $\{i_1,
i_2,i_3,i_4\} \subset \{1, \dots,5\} $ such that $\{i_1, i_2,i_3,i_4\}\in H_8$.
Let $a= \{1, \dots,8\}\setminus \{i_1, i_2,i_3,i_4\}$. Then $|a\cap \{1,
\dots,5\}|=1$ and $(h-1)\a_a$ is a root.

\textbf{Case 3:} $L^-(h)\cong D_4$. Then $k=4$ and $\{1,2,3,4\}\in H_8$. Since
$H_8$ is a self dual code, for any $a\in H_8$, $|a\cap\{1,2,3,4\}|$ is even.
Hence, for any $\a\in E_8\setminus A$,  $(h-1)\a$ is either $0$ or has $2$ or
$4$ non-zero projections to the $A_1$'s. Thus, by Lemma \refpp{2non0} (a), $h$
is rootless.

\textbf{Case 4:} $L^-(h)\cong D_6, E_7$ or $E_8$.  Then $k\geq 6$. Since the minimal weight of $H_8$ is $4$,
we have $|a\cap \{1, \dots, k\}|\geq 2$ for any nonzero element $a\in H_8$. Hence, $h$ is rootless by Lemma  \refpp{2non0} (a).
\eop

\medskip

\noindent {\bf Order 4}

\begin{nota}\labtt{order4j}
Let $h$ be a rootless element of order 4 and set $J:=\ker(h^2+1)$.
\end{nota}

Then $J$
has even rank and $h^2$ is also rootless. Since $det(h^2)=1$,  \refpp{order2} implies that  $J\cong D_4, D_6$, or $E_8$.

\begin{lem}\labtt{order4h}
Let $h\in O(D_{2n})$ be an element of order 4 and $h^2=-1$. Then there exists
an orthogonal set of roots $\{\a_1, \dots,\a_{2n}\}\subset D_{2n}$ such that $h
(\a_{2i-1}) =\a_{2i}$ and $h(\a_{2i})= -\a_{2i-1}$ for all $i=1,\dots, n$.
\end{lem}

\pf We shall use the standard model for $D_{2n}$, i.e.,
\[
D_{2n}=\{ \sum_{i=1}^{2n}  x_i e_i  \mid x_1+\cdots+x_{2n}\equiv 0 \mod 2\},
\]
where $\{e_1, \dots, e_{2n}\}$ is the standard basis of $\ZZ^{2n}$.

Then up to conjugacy in $O(D_{2n})$, we may assume that $h=DP$, where $P$ is a
matrix associated to a permutation $\sigma\in Sym_{2n}$ and $D$ is a diagonal
matrix with diagonal entries $1$ or $-1$. Note that $$P=\sum_{i=1}^{2n}
E_{\sigma i, i},$$ where $E_{i,j}$ is a matrix whose $(i,j)$-entry  is $1$ and
all other entries are $0$.

Let $\epsilon_1, \dots, \epsilon_{2n}$ be the diagonal entries of $D$. Then
\[
DPD= \sum_{i=1}^{2n} \epsilon_{\sigma i}\epsilon_{i} E_{\sigma i, i}
\] and
\[
(DP)(DP) =(DPD)P = \sum_{1\leq i,j \leq 2n} \epsilon_i \epsilon_{\sigma_i}  \delta_{i,\sigma_j} E_{\sigma i, j}.
\]

By $h^2=-1$, we have $(DP)(DP)=(DPD)P =- I$. This implies  $\sigma^2 =1$ and
$\epsilon_{\sigma_i} \epsilon_{i}=-1$. Therefore, by rearranging the indices if
necessary, the matrix of $h$ with respect to the standard basis is given by
\[
\begin{pmatrix}
0& 1 &0 &0 & \dots & 0& 0 \\
-1&0 &0&0 & \dots & 0& 0  \\
0&0&0&1 & \dots & 0& 0 \\
0&0&-1 &0 & \dots & 0& 0\\
\vdots &\vdots&\vdots&\vdots& \ddots& \vdots&\vdots \\
0&0&0&0 & \dots & 0& 1 \\
0&0&0 &0 & \dots & -1& 0
\end{pmatrix}.
\]

Now define $\a_{2i-1} = e_{2i-1} - e_{2i}$ and $\a_{2i}= e_{2i-1} +e_{2i}$ for $i=1,
\dots n$. Then $\{\a_1, \a_2, \dots, \a_{2n-1}, \a_{2n}\}$ satisfies the required
properties. \eop

\bigskip

We now treat the order 4 case according to the three types of $J$
\refpp{order4j}.
\begin{nota}\labtt{order4f}
Let $F:=ann_{E_8}(J)$. Note that $h^2$ acts trivially on $F$.
\end{nota}
\smallskip

\textbf{Case 1:} $J\cong E_8$. Then $h$ is fixed point free and $h^2$ acts as
$-1$ on $E_8$.   Such elements form one conjugacy class \refpp{order4h}.

\medskip

\textbf{Case 2:}  $J\cong D_6$. Then $F\cong A_1\perp A_1$. Then by Lemma
\ref{order4h}, there exists $\{\a_1, \a_2, \dots, \a_6\}\subset J$ such that
$h(\a_{2i-1})= \a_{2i}$, $ h(\a_{2i})=- \a_{2i-1}$ for $i=1,2,3$ and
\[
J=span_\ZZ \{ \a_1, \dots, \a_6, \frac{1}2(\a_1+\a_2+\a_3+\a_4),
\frac{1}2(\a_3+\a_4 + \a_5+\a_6)\}.
\]

Let $\{\a_7, \a_8\}$ be a basis of $F$. Then we may also arrange indexing so
that
\[
E_8= span_\ZZ \left \{ {\a_1, \dots, \a_8, \frac{1}2(\a_1+\a_2+\a_3+\a_4),
 \frac{1}2(\a_3+\a_4+\a_5+\a_6),  } \atop { \frac{1}2( \a_5+\a_6+\a_7+\a_8), \frac{1}2(\a_1+\a_3 + \a_5+\a_7)} \right\}.
 \]

 Next we shall study the action of $h$ on $F$.

 \begin{lem}  In above notation,
$ h(\a_7)\in span_\ZZ \{\a_8\}$.
 \end{lem}

\pf
 Suppose $h(\a_7)\notin span_\ZZ \{\a_8\}$. Then $h(\a_7) =\pm \a_7$  and $h(\a_8)=\pm \a_8$.

 In this case, we have
 \[
 (h-1) \frac{1}2(\a_1+\a_3+\a_5+\a_7)  = \frac{1}2( -\a_1+\a_2-\a_3+\a_4-\a_5+\a_6- \a_7+\epsilon \a_7) , \quad \epsilon =\pm 1,
 \]
 which has norm 3 or 5. It is a contradiction since $E_8$ is even.
\eop

By the lemma above, we may assume $h(\a_7)=\a_8$ and $h(\a_8)=\a_7$ (by
replacing $\a_8$ by $-\a_8$ if necessary).  Then
\[
(h-1) \frac{1}2 (\a_5+\a_6+\a_7+\a_8) = \a_5,
\]
which is a root. Thus, $h$ is  not rootless.

\medskip

\textbf{Case 3:}  $J\cong D_4$ and $F\cong D_4$.  This will lead to two cases
for $h$.

\begin{nota}\label{ai}
Let $\{\a_1, \a_2, \a_3,\a_4\}\subset J$  such that $h(\a_1)=\a_2,
h(\a_2)=-\a_1, h(\a_3)=\a_4, h(\a_4)=-\a_3$ (cf. Lemma \refpp{order4h}).

Let $\{\a_5, \a_6, \a_7, \a_8\}\subset F$ such that $( \a_i,\a_j)
=2\delta_{i,j}$.

We may reindex to assume
\[
E_8= span_\ZZ \left \{ {\a_1, \dots,, \a_8,  \frac{1}2(\a_1+\a_2+\a_3+\a_4),
 \frac{1}2(\a_3+\a_4+\a_5+\a_6),  } \atop { \frac{1}2( \a_5+\a_6+\a_7+\a_8), \frac{1}2(\a_1+\a_3 + \a_5+\a_7)} \right\}.
 \]
\end{nota}

\smallskip

 \begin{lem}\label{lem:1}
If $h$ is rootless, then $h(\a_i)=\pm \a_i$ for all $i=5, \dots,8$.
 \end{lem}

\pf Suppose $h(\a_k)= \epsilon \a_\ell$ for some $\epsilon =\pm 1$, $k\neq \ell$
and  $k,\ell\in \{5,6,7,8\}$. Then $h(\epsilon \a_\ell)=h^2(\a_k)=\a_k$ since
$\a_k\in F$.

Take $i,j\in \{1,2,3,4\}$ with $i<j$ such that
\[
\frac{1}2(\a_i+\a_j+\a_k +\epsilon \a_\ell)\in E_8.
\]
Then
\[
\begin{split}
&\  h\left( \frac{1}2(\a_i+\a_j+\a_k +\epsilon \a_\ell)\right)\\
= &
\begin{cases}
\frac{1}2(\a_i-\a_j+\a_k +\epsilon \a_\ell) & \text{ if }  h(\a_i)\in span_\ZZ(\a_j), \\
\frac{1}2(\pm \a_{i'} \pm \a_{j'} + \a_k +\epsilon \a_\ell)& \text{ if }
h(\a_i)\notin span_\ZZ(\a_j),
\end{cases}
\end{split}
\]
where $\{\a_i, \a_j, \a_i', \a_j'\} =\{\a_1, \a_2, \a_3,\a_4\}$.

In either case, $(h-1) \frac{1}{2} (\a_i+\a_j+\a_k +\epsilon \a_\ell)$ is a root. \eop

\begin{lem}
Let $Y$ be the fixed point sublattice of $h$ on $F$. Then $rank\, Y\leq 1$.
\end{lem}

\pf Suppose  $rank\, Y\geq 2$. Then by the previous lemma, $h$ fixes $\a_k$ and
$\a_\ell$ for some   $k\neq \ell$ and  $k,\ell\in \{5,6,7,8\}$. Take $i,j\in
\{1,2,3,4\}$ with $i<j$ such that
\[
\frac{1}2(\a_i+\a_j+\a_k + \a_\ell)\in E_8.
\]
Then by the same argument as in Lemma \ref{lem:1}, $
(h-1)\frac{1}2(\a_i+\a_j+\a_k +  \a_\ell)$ is a root. \eop

Since $h(\a_i)=\pm \a_i$ for $i=5, \dots,8$ and $rank\, Y\leq 1$, $\{\a_5,\a_6\}$
or $\{\a_7, \a_8\}$ is contained in the $(-1)$-eigenspace of $h$.

By reindexing, we may assume $\a_5, \a_6$ are in the $(-1)$-eigenspace of $h$.
Define
\[
\begin{split}
\b_1&:=\frac{1}2 (\a_1+\a_2+\a_5+\a_6),\\
\b_1'&:= \frac{1}2 (\a_3+\a_4+\a_5-\a_6) = \frac{1}2 (\a_3+\a_4+\a_5+\a_6) -\a_6.
\end{split}
\]
Then by our convention \refpp{ai}, $\b_1$ and $\b_1'$ are in $E_8$. Let
\[
\b_2:=h(\b_1)=\frac{1}2 (-\a_1+\a_2-\a_5-\a_6),\quad  \b_3:=h^2(\b_1)=\frac{1}2
(-\a_1-\a_2+\a_5+\a_6),
\]
\[
\b_2':=h(\b_1')=\frac{1}2 (-\a_3+\a_4-\a_5+\a_6), \quad \b_3':=h^2(\b_1')=\frac{1}2 (-\a_3-\a_4+\a_5-\a_6).
\]
Then $\b_2, \b_3, \b_2',\b_3'$ are also in $E_8$ since $h\in O(E_8)$.

Let $A:=span\{\b_1, \b_2,\b_3\}$ and $A':= span\{\b_1', \b_2',\b_3'\}$. Then
$A\cong A'\cong A_3$ and $(A, A')=0$. By identifying $A, A'$ with $A_3$,
$h|_A$ and $h|_{A'}$ are identified with $h_{A_3}$.

Let $X:= ann_F( span\{\a_5, \a_6\})$. Then $X\cong A_1\oplus A_1$ and
$Y=Fix_F(h) < X$. Note also that $(X,A)=(X,A')=0$.

If $Y=0$, then $h|_{X}=-id_X$. If $Y\cong A_1$, then $h$ acts trivially on $Y$
and acts as $-1$ on $X':=ann_{X}(Y)\cong A_1$. Thus, $h$ may be identified with
\[
\begin{cases}
h_{A_3}\oplus h_{A_3}\oplus h_{A_1}\oplus h_{A_1} & \text{  if }  Y=Fix(h)=0,\\
h_{A_3}\oplus h_{A_3}\oplus h_{A_1}\oplus id_{A_1} & \text{  if }  Y=Fix(h)\cong A_1.
\end{cases}
\]

\medskip

Let $Q\cong A_3\oplus A_3\oplus A_1\oplus A_1$ be a sublattice of $E_8$. Then
$|E_8/Q|=8$ and any element in $E_8\setminus Q$ has non-zero projections to at
least three $A_i$'s. If $h=h_{A_3}\oplus h_{A_3}\oplus h_{A_1}\oplus h_{A_1}$
or $h_{A_3}\oplus h_{A_3}\oplus h_{A_1}\oplus id_{A_1}$, then $(h-1)x$, $x\in
E_8\setminus Q$, has at least two non-zero projections to the $A_i$'s.
Therefore, they are rootless by \refpp{2non0}.

As a summary, we have
\begin{thm}\labtt{order4}
Let $h$ be a rootless element of order 4. Then $J=\ker(h^2+1)\cong D_4$ or
$E_8$.

(1) If $J\cong D_4$, then $h$ conjugate to $h_{A_3}\oplus h_{A_3}\oplus
h_{A_1}\oplus h_{A_1}$ or $h_{A_3}\oplus h_{A_3}\oplus h_{A_1}\oplus
id_{A_1}$.

(2) If $J\cong E_8$, then  $h$ is fixed point free and $h^2$ acts as $-1$ on
$E_8$.   Such elements form one conjugacy class.
\end{thm}

\noindent {\bf Order 8}

\begin{thm}\labtt{rootless8}
There is no rootless element of order $8$.
\end{thm}

\pf Suppose $h$ is a rootless element of order 8. Then $g=h^2$ is a rootless
element of order $4$. By the analysis of order 4 elements, $Ker (g^2+1)\cong
D_4$ or $E_8$ (cf. Theorem \ref{order4}).

In either case, there exists  a $D_4$ sublattice of $E_8$ which $h$ acts (cf.
Lemma \ref{order4h}).

Recall that $O(D_4)$ has the shape $(2^3{:}Sym_4).Sym_3$ (see (4.3.12) in \cite{gal})
Since $h$ has order $8$, $h$ acts on $D_4$ as a product of
a 4-cycle in $Sym_4$ and an outer involution with respect to the standard model of $D_4$.
Therefore, there exists $\{\a_1,\a_2, \a_3,\a_4\}$ such that
$(\a_i,\a_j)=2\delta_{i,j}$ for $i=1,2,3,4$ and
\[
h(\a_1)=\a_2,\ h(\a_2)=\a_3,\ h(\a_3)=\a_4,\ h(\a_4)=-\a_1.
\]
However,
\[
(h-1) \frac{1}2( \a_1+\a_2+\a_3+\a_4) = - \a_1,
\]
which is root, a contradiction. \eop

\subsection{Rootless elements of composite orders }

\noindent
{\bf Order 6}

Let $h$ be a rootless element of order $6$. Let $g:=h^2$ and $t
:= h^3$. Then, $g$ has order $3$ and $t$ has order $2$.

Let $L^+(t)$ and $L^-(t)$ be the $(+1)$ and  $(-1)$-eigenlattice of $t$ on
$E_8$.

\begin{lem}
If $h$ is rootless of order 6, then $L^+(t)\cong D_4$.
\end{lem}

\pf First, we note that $g=h^2$ acts on both $L^+(t)$ and $L^-(t)$.

By the order 2 analysis,  $L^+(t)\cong 0, A_1, A_1^2$, or $D_4$.

\noindent \textbf{Case 1:} $L^+(t)=0$  and thus $t$ acts as $-1$ on $E_8$.
Therefore,
\[
h=t g^2= -g^2,
\]
By the order 3 analysis,  we may identify $g^2$ with either $h_{A_2}^{\oplus
4}$ or $h_{A_2}^{\oplus 3}\oplus id_{A_2}$.

In either case, let $\hat{\gamma}=(\gamma_1, \gamma_2, \gamma_3,0)$ be a
root in $E_8$, where $\gamma_1, \gamma_2, \gamma_3\in A_2^*$ and have
norm $2/3$.  Since $1+h_{A_2}+h_{A_2}^2=0$ on $A_2^*$,
$(-h_{A_2}-1)\gamma_i = h^2_{A_2} \gamma_i$ also has norm $2/3$ for
$i=1,2,3$. Therefore,
\[
(h-1)(\hat{\gamma}) = ((-h_{A_2}-1)\gamma_1, (-h_{A_2}-1)\gamma_2, (-h_{A_2}-1)\gamma_3,0)
\]
has norm $2$ and is a root.

\medskip

\noindent \textbf{Case 2:}  $L^+(t)\cong A_1$. Then $g$ acts trivially on
$L^+(t)$. Thus, $Fix(g)\neq 0$ and hence $Fix(g)\cong A_2$ and $g^2$ may be
identified with $h_{A_2}^{\oplus 3}\oplus id_{A_2}$ by Theorem \ref{rootless3}.
Note that $L^+(t)< Fix(g)$. Therefore, $ann_{E_8}(Fix(g)) <  ann_{E_8}(L^+(t))
= L^-(t)$ and we have
\[
h|_{ann_{E_8}(Fix(g))} =- g^2.
\]
Let $\hat{\gamma}= (\gamma_1, \gamma_2, \gamma_3,0)$ be a root  in $
ann_{E_8}(Fix(g)) \cong E_6$, where $\gamma_1, \gamma_2, \gamma_3\in
A_2^*$ have norm $2/3$ . Then, as in Case 1,
\[
(h-1)\hat{\gamma}= ((-h_{A_2}-1)\gamma_1, (-h_{A_2}-1)\gamma_2, (-h_{A_2}-1)\gamma_3,0)
\]
is a root.

\medskip

\noindent \textbf{Case 3:}  $L^+(t)\cong A_1\oplus A_1$. Then $g$ acts trivially
on $L^+(t)$ since $O(A_1\oplus A_1)$ has no elements of order 3.  This is
impossible since $Fix(g)\cong A_2$ does not contain a sublattice of type
$A_1+A_1$.

Therefore, the only possible case is  $L^+(t)\cong D_4$. \eop

\medskip

Since $L^+(t)\cong D_4$, we also have $L^-(t)=ann_{E_8}(L^+(t))\cong
ann_{E_8}(D_4)\cong D_4$ \cite[(5.3.1)]{gal}. Note that $g$ acts on both
$L^+(t)$ and $L^-(t)$.

\begin{lem}\labtt{evenrank}
Let $Fix_{L^\pm (t)}(g)$ be the fixed points of $g$ on $L^\pm (t)$. Then the rank
of $Fix_{L^\pm (t)}(g)$ is even.
\end{lem}

\pf Note that the minimal polynomial of $g$ on $ann_{L^\pm(t)}( Fix_{L^\pm
(t)}(g))$ is $x^2+x+1$, which is irreducible. Thus $rank(ann_{L^\pm(t)}(
Fix_{L^\pm (t)}(g)))$ is even and so is  $rank(Fix_{L^\pm (t)}(g))$. \eop

\begin{lem}\labtt{nonzero}
 We use the same notation as in \refpp{evenrank}. Then $Fix_{L^- (t)}(g)\neq 0$.
\end{lem}

\pf Suppose $g$ is fixed point free on $L^-(t)$. Then $span\{\a, g\a\}\cong A_2$
for any root $\a\in L^-(t)$. Now choose a root $\a\in L^-(t)$ and define
$A:=span\{\a, g\a\}$.

Let $B:=ann_{L^-(t)}(A)$. Then $B\cong \sqrt{2}A_2$.  Thus,  we obtain a
sublattice $A\oplus B \cong A_2\oplus\sqrt{2}A_2$ in $L^-(t)$  and  $g$ acts
fixed point freely on the indecomposable direct summands.

By the previous lemma, $Fix_{L^+(t)}(g)$ has even rank and hence
$Fix_{L^+(t)}(g)\cong A_2$ or $0$.
We shall first obtain information in these two cases, then finally a contradiction to prove this lemma.

\textbf{Case 1:} $X:=Fix_{L^+(t)}(g)\cong A_2$. Then $C:= ann_{L^+(t)}(X)\cong
\sqrt{2}A_2$ and $g$ acts fixed point freely on $C$. Thus, we obtain a sublattice
$$ X \oplus A\oplus B\oplus C \cong  A_2\oplus A_2\oplus  \sqrt{2}A_2\oplus \sqrt{2}A_2 $$
in $E_8$  such that $g$ acts on each indecomposable summand and is fixed
point free on $B$ and $C$.

Notice that $B\oplus C <  ann_{E_8}(X\oplus A) \cong A_2\oplus A_2$ and
$$|ann_{E_8}(X\oplus A) /(B\oplus C)|= 2^2.$$

Since $ann_{E_8}(X\oplus A) \cong A_2\oplus A_2$ has roots, there exist $\b\in
B$ and $\gamma\in C$ with $(\b,\b)=(\gamma,\gamma)=4$ such that
$\frac{1}2(\b+\gamma)$ is a root in $ann_{E_8}(X\oplus A)$.  Then
 we also have  $\frac{1}2(g\b+g\gamma)\in ann_{E_8}(X\oplus A)$.
 Recall that the 2-part of $\dg{\sqrt{2}A_2}=(\sqrt{2}A_2)^*/ \sqrt{2}A_2$
 is generated by the elements of the
form $\frac{1}2 \d +\sqrt{2}A_2$ for $\d\in \sqrt{2}A_2$ with $(\d,\d)=4$.

 By comparing
the determinants, we have
\[
ann_{E_8}(X\oplus A)
= span\{ B\oplus C, \frac{1}2(\b+\gamma), \frac{1}2(g\b+g\gamma)\} \cong A_2\oplus A_2.
\]

Let $A^+= span\{\frac{1}2(\b+\gamma),\frac{1}2(g\b+g\gamma)\}$ and $A^-=
span\{\frac{1}2(-\b+\gamma), \frac{1}2(-g\b+g\gamma)\}$. Then $A^+$ and $A^-$
are sublattices of $ann_{E_8}(X\oplus A)$. Since $g$ satisfies $x^2+x+1=0$ on
$ann_L(X)$, we have $(v,gv)=-\half (v,v)$ for all $v\in ann_L(X)$. It follows
that $A^+\cong A^-\cong A_2$ and $(A^+, A^-)=0$. Moreover, $g$ stabilizes each
of $A^+$ and $A^-$.

Note that $t$ commutes with $g$ and $h=tg^2$. Since $X, C<  L^+(t)$ and $A,B <
L^-(t)$ , we have
\[
h|_X=id_X,\quad  h|_{A} =-g^2|_{A},
\]
\[
h(\frac{1}2(\b+\gamma))= \frac{1}2(-g^2\b + g^2\gamma),
\quad h(\frac{1}2(-\b+\gamma))= \frac{1}2(g^2\b + g^2\gamma).
\]
Thus we have $h(A^+)=A^-$ and $h(A^-)=A^+$.  Note that
\[
t(\frac{1}2(\b+\gamma)) = \frac{1}2(-\b+\gamma) \text{\ and \ }
t(\frac{1}2(g\b+g\gamma)) = \frac{1}2(-g\b+g\gamma).
\]
Therefore, $h$ acts on $A^+\oplus A^-$ and $t$ interchanges  $A^+$ and $A^-$.

By identifying $X\oplus A\oplus A^+\oplus A^-$ with $A_2^4$ and $g^2$ with
$h_{A_2}$ on $A, A^+$ and $A^-$, $h$ is conjugate to $\sigma\tau$, where
\[
\sigma =  id_{A_2} \oplus (-h_{A_2}) \oplus h_{A_2}\oplus h_{A_2}
\]
and $\tau$ performs a transposition on the 3rd and 4th copies of $A_2$ and is
the identity on the first two summands. 

\textbf{Case 2:} $Fix_{L^+(t)}(g)=0$. Then $g$ acts fixed point freely on
$Fix_{L^+(t)}(g)$. Let $\a\in L^+(t)$ be a root. Then $X':= span \{ \a,
g\a\}\cong A_2$. Let $C':= ann_{L^+(t)}(X)$. Then $C'\cong \sqrt{2}A_2$ and we
obtain a sublattice $X'\oplus A\oplus B\oplus C'\cong A_2\oplus A_2\oplus
\sqrt{2} A_2\oplus \sqrt{2}A_2$ in $E_8$ such that $g$ acts fixed point freely
on $X'$, $A$, $B$ and $C'$. Then by an argument as in case 1, one can show that
$h$ is conjugate to $\sigma'\tau$, where
\[
\sigma' =  h_{A_2} \oplus (-h_{A_2}) \oplus h_{A_2}\oplus h_{A_2}
\]
and $\tau$ is a transposition on the 3rd and 4th copies of $A_2$.

We now get a contradiction to both Case 1 and Case
 2.
We take a sublattice $A_2^4$ of $E_8$ so that $g$ preserves each summand and
$h$ has the form $\s \t$, $\s' \t$, as described in the two cases. Let $\eta:=
\frac{1}3(0,a,b,c)$ be a root in $E_8$ where $a,b,c\in A_2$ have norm 6. Then,
$h\eta = \frac{1}3(0,-h_{A_2}a, h_{A_2}c, h_{A_2}b)$ and
\[
(\eta, h\eta)= \frac{1}9\left( 3 +(b,h_{A_2}c)+ (c,h_{A_2}b)\right)= \frac{1}9 (3- (b,c))
\]
since $(1+h_{A_2}+ h_{A_2}^{2})b=0$ and $(c,h_{A_2}b)= (b,h_{A_2}^{2}c)$.

Since $\eta$ is a root, $(\eta, h\eta)=0, \pm 1$ or $\pm 2$. Thus, we have
$(b,c)=-6$ or $3$ because $|(b,c)|\leq 6$  and $\frac{1}9 (3- (b,c))\in \ZZ$.
It implies $c=-b $ or $-h_{A_2}^i b$ for $i=1,2$.

Since $h_{A_2}$ stabilizes all cosets of $A_2$ in $A_2^*$, we also have $
\frac{1}3(0,a,b,h_{A_2}^i c)\in E_8$ for all $i=1,2$. Thus, by replacing $c$ by
$h_{A_2}^i c$ if necessary, we may assume $c=-b$. Then
\[
(h-1)\eta = -\frac{1}3 (0,(h_{A_2}+1)a, (h_{A_2}+1)b, (h_{A_2}+1)c).
\]
Recall that $(h_{A_2}\a, \a)= -\frac{1}2 (\a,\a)$ for $\a=a,b,c$ (cf.
\cite[Lemma 3.2]{GLEE8}) Thus, $ (h_{A_2}+1)a, (h_{A_2}+1)b$ and $(h_{A_2}+1)c$
have norm 6 and $(h-1)\eta$ is a root. This final contradiction proves that
$Fix_{L^- (t)}(g)\neq 0$. \eop

\begin{lem}\labtt{fixg}
We use the same notation as in \refpp{evenrank} and \refpp{nonzero}. Then
$Fix_{L^-(t)}(g)\cong A_2$ and $g$ acts fixed point freely on $L^+(t)$.
\end{lem}
\pf We first note that $Fix_L(g)\cong A_2$ or $0$ (see \refpp{rootless3}).
Since $Fix_{L^-(t)}(g)\neq 0$ and has even rank,  we have $Fix_{L^-(t)}(g)\cong
A_2$ and $Fix_{L^+(t)}(g)=0$. \eop

\medskip

By the same argument as in Lemma \refpp{nonzero}, we have the following.
\begin{lem}\label{st}
Let $h$ be a rootless element of order $6$. Then $h$ is conjugate to
$\sigma\tau = \tau \sigma$, where
$\sigma = (-id_{A_2}) \oplus h_{A_2}  \oplus h_{A_2}\oplus
h_{A_2}$ and $\tau$ is an involution which interchanges
the 3rd and 4th copies of $A_2$.
\end{lem}
\pf Let $P:= Fix_{L^-(t)}(g)\cong A_2$ and $R:= ann_{L^-(t)}(P) (\cong
\sqrt{2}A_2)$. Take a root $\a\in L^+(t)$. Then $Q:= span\{ \a, g\a\} \cong A_2$
since $g$ acts fixed point freely on $L^+(t)$.  Also,  $S:= ann_{L^+(t)}(Q)\cong
\sqrt{2}A_2$. Thus we obtain a sublattice $P\oplus Q\oplus R\oplus S\cong
A_2\oplus A_2\oplus \sqrt{2}A_2 \oplus \sqrt{2}A_2$ in $E_8$ such that $g$
acts trivially on $P$ and fixed point freely on $Q,R$ and $S$. Again, we have
$R\oplus S < ann_{E_8} (P\oplus Q)\cong A_2\oplus A_2$. Thus, by the same
argument as in Lemma \refpp{nonzero}, one can show that $h$ is conjugate to
$\sigma\tau$, where $\sigma = (-id_{A_2}) \oplus h_{A_2} \oplus h_{A_2}\oplus
h_{A_2}$ and $\tau$ is an involution which interchanges the 3rd and 4th copies
of $A_2$.\eop

\medskip

Let $\sigma$ and $\tau$ be as in Lemma \ref{st} and assume $h=\sigma\tau$. Then
we determine a sublattice $(A_2)^4$ in $E_8$.

Let  $\eta:= \frac{1}3(\b, 0 , \gamma, \gamma')\in (A_2^*)^4$ be a root in
$E_8$, where $\b$, $\gamma$ and $\gamma'$ have norm $6$. Then $h(\eta)=
\frac{1}3(-\b, 0, h_{A_2} \gamma',h_{A_2}\gamma)$ and
\[
(\eta, h\eta)= \frac{1}9\left( (\b,-\b) +(\gamma,h_{A_2}\gamma')+ (\gamma' ,h_{A_2}\gamma)\right)
= \frac{1}9 (-6 - (\gamma, \gamma')).
\]
Since $(\eta, h\eta)=0, \pm 1$ or $\pm2$, we have $(\gamma, \gamma')=3$ or $-6$
and hence $\gamma' = -h_{A_2}^i \gamma$ for $i=0,1,2$. Without loss, we may
assume $\gamma'= - h_{A_2}\gamma$ since $h_{A_2}$ stabilizes all cosets of
$A_2$ in $A_2^*$.

Then, we have $\eta=\frac{1}3(\b, 0, \gamma, -h_{A_2} \gamma)$ and
\[
\begin{split}
h\eta &=\frac{1}3(-\b, 0, -h_{A_2}^2 \gamma, h_{A_2} \gamma),\\
h^2\eta&=
\frac{1}3(\b, 0, h_{A_2}^2 \gamma, -\gamma), \\
h^3\eta &= \frac{1}3(-\b,0,  -h_{A_2} \gamma,\gamma),\\
h^4\eta &= \frac{1}3(\b,0 , h_{A_2} \gamma, - h_{A_2}^2\gamma).
\end{split}
\]
Thus  we have $( h\eta, \eta) = ( h^{-1}\eta , \eta) =-1$, $( h^2\eta, \eta) = (
h^{-2} \eta , \eta) =0$ and $( h^3\eta, \eta)=0$. It implies that $ A=span\{ h^i
\eta \mid i=0, \dots, 5\} \cong A_5 $ and $\{\eta, h\eta, h^2\eta, h^3\eta,
h^4\eta\}$ is a fundamental set of simple roots. By identifying $A$ with $A_5$,
we may identify $h|_A$ with $h_{A_5}$.

Let $B$ be the second summand isometric to $A_2$ and $C:= ann_{L^-(h)}(\b)$.
Then $C\cong A_1$ and $h$ acts as $-1$ on $C$. Thus we have  a rank 8
sublattice $A\oplus B\oplus C$ in $E_8$ such that $A\cong A_5$, $B\cong A_2$,
$C\cong A_1$. Moreover, we may identify  $h|_{A}$ with $h_{A_5}$, $h|_B$ with
$h_{A_2}$ and $h|_{C}=-id_C$. The following theorem now follows.

\begin{thm}\labtt{order6}
Let $h$ be a rootless element of order 6. Then $h$ is conjugate to
$h_{A_5}\oplus h_{A_2}\oplus h_{A_1}$.
\end{thm}

\medskip

\noindent\textbf{Other composite orders}

\begin{thm}\labtt{no12}
There is no rootless element of order 12.
\end{thm}
\pf Let $h$ be a rootless element of order 12. Then $g=h^4$ has order 3,
$f=h^3$ has order 4 and both are rootless. By the analysis of rootless order 6
elements, we have $Fix_{L^-(f^2)}(g)\cong A_2$ (see \refpp{fixg}). Since $f$
commutes with $g$, $f$ also acts on $Fix_{L^-(f^2)}(g)$. For any root $\a\in
L^-(f^2)$, we have
\[
(f\a,\a)=(f^2\a, f\a)=-(\a, f\a).
\]
Hence $(f\a, \a)=0$ and $span\{\a, f\a\}\cong A_1\oplus A_1$. Since $A_2$ does
not contain any sublattice isometric to $A_1\oplus A_1$, $f$ cannot stabilize
any $A_2$-sublattice in $L^-(f^2)$, which is a contradiction. \eop

\begin{lem} If $h\in O(L)$ is rootless, $|h|$ is not 10 or 15.
\end{lem}
\pf Let $h$ be rootless and have order 10 or 15.  We use the notations in \refpp{ppart}.   Since $h_5$ is fixed point free, if
$q$ is the other prime dividing $|h|$, the $q$-part has eigenvalue 1.   This
means if $q=3$, then $h_3$ has rank 2 fixed point sublattice, which is impossible
since $h_5$ does not leave invariant a rank 2 sublattice. Now suppose that
$q=2$. Since the fixed point sublattice $F$ of $h_2$ is nonzero and is
$h$-invariant, $rank(F)=4$.  However, no rank 4 RSSD sublattice of $L$ has an
automorphism of order 5, contradiction. \eop

\section{How the surviving cases give all rootless $EE_8$ pairs}\labttr{sec:2.1}

Each of the 11 lattices from the main result of \cite{GLEE8} has the form
$M+N$, where $M\cong N \cong EE_8$ and is denoted by some notation $DIH_{2k}(d,
\cdots )$, where $d$ is the rank and $2k=|\la t_M, t_n \ra|$. Their structures
are summarized in Table 1. We shall prove that each of the 11 cases occurs as
some SDC-lattice $L(E_8,h)$ by using the rootless $h$, which we classified in
preceding sections.

We exclude the case $h=1$, which is indeed rootless, but for which $M=N=L$.

\begin{center}
Table 1: {\bf Integral rootless lattices which are sums of $EE_8$s}\\
\smallskip
\begin{tabular}{|c|c|l|c|c|}
\hline
 Name & $\<t_M,t_N\>$ & Isometry type of $L$ (contains)& $\dg L$ & In Leech? \cr
 \hline \hline
 $\dih{4}{12}$ & $Dih_4$ & $\geq  DD_4^{\perp 3}$ & $1^4 2^6 4^2$ & Yes\cr
 \hline
 $\dih{4}{14}$ & $Dih_4$ & $\geq AA_1^{\perp 2} \perp DD_6^{\perp 2}$ & $1^4 2^8 4^2$  &Yes\cr
 \hline
 $\dih{4}{15}$ & $Dih_4$ & $\geq AA_1\perp EE_7^{\perp 2}$ & $1^22^{14}  $  & No\cr
 \hline $\dih{4}{16}$ & $Dih_4$ & $\cong EE_8\perp EE_8$ & $2^{16}$ &   Yes\cr
 \hline $\dih{6}{14}$ & $Dih_6$ & $\geq AA_2\perp A_2\otimes E_6 $ & $1^7 3^3 6^2$ &Yes\cr
 \hline
 $\dih{6}{16}$ & $Dih_6$ & $\cong A_2\otimes E_8 $ & $1^8 3^8 $ &Yes\cr
 \hline
 $\dih{8}{15}$ & $Dih_8$ & $ \geq AA_1^{\perp 7}\perp EE_8$     & $1^{10}4^5 $&   Yes\cr
 \hline
 $\dih{8}{16, DD_4}$ & $Dih_8$ & $\geq DD_4^{\perp 2} \perp  EE_8$  & $1^8 2^4 4^4$  &Yes\cr
 \hline
 $\dih{8}{16, 0}$ & $Dih_8$ & $\cong BW_{16} $ & $1^8 2^8$ &Yes\cr
 \hline $\dih{10}{16}$ & $Dih_{10}$ & $\geq A_4\otimes A_4 $ & $1^{12}5^4$  &Yes\cr
 \hline $\dih{12}{16}$ & $Dih_{12}$ & $\geq AA_2\perp AA_2 \perp A_2\otimes E_6 $ & $1^{12}6^4 $  & Yes\cr
 \hline
\end{tabular}
\\ \smallskip
$X^{\perp n}$ denotes the orthogonal sum of $n$ copies of the
lattice $X$.
\end{center}

There are 11 rootless nonidentity conjugacy classes.  If we form the associated 11 SDC lattices, it suffices to argue that they give 11 distinct $EE_8$-pairs.  Notice that the dihedral group $\la t_M, t_N\ra$ has order $2|h|$ \refpp{sdc3}.

We now prove the bijection by use of Table 2. In column 1, we list the
possibilities for rootless $h$.  Columns 2 and 3 are consequences of our
classification of rootless elements of $O(E_8)$. Our intended correspondence is
expressed in column 4, which we shall now justify.

\begin{center}
Table 2: {\bf Rootless classes in $O(E_8)$}\\
\begin{tabular}{|c|c|c|c|}
\hline
Notation for $h$ & Order of $\la t_M, t_N\ra$ &  $rank(M+N)$ &
Lattice name in \cite{GLEE8} \cr
\hline
$h_{A_1}^8$ & $4$ & $16 $ & $ \dih{4}{16}$ \cr \hline
$h_{A_1}^7\oplus id_{A_1}$ & $4$ & $15 $ & $\dih{4}{15}$ \cr \hline
$h_{A_1}^6\oplus id_{A_1^2}$ & $4$ & $14 $ & $\dih{4}{14}$ \cr \hline
$h_{A_1}^4\oplus id_{A_1}$ & $4$ & $12 $ & $ \dih{4}{12}$ \cr \hline
$h_{A_2}^4$ & $6$ & $16 $ & $ \dih{6}{16}$ \cr \hline
$h_{A_2}^3\oplus id_{A_2}$ & $6$ & $14 $ & $ \dih{6}{14}$ \cr \hline
$h_{A_3}^2\oplus h_{A_1}^2$ & $8$ & $16 $ & $\dih{8}{16, DD_4}$ \cr \hline
$h_{A_3}^2\oplus h_{A_1}\oplus id_{A_1}$ & $8$ & $15 $ & $\dih{8}{15}$ \cr \hline
$h^2=-1$ & $8$ & $16 $ & $\dih{8}{16,0}$ \cr \hline
$h_{A_4}^2$ & $10$ & $16 $ & $\dih{10}{16}$ \cr \hline
$h_{A_5}\oplus h_{A_2} \oplus  h_{A_1}$ & $12$ & $16 $ & $ \dih{12}{16}$ \cr
\hline
\end{tabular}
\end{center}

We observe that two lattices which occur for different  entries in column 1 of
Table 1 are distinguished by the orders of the dihedral groups and their ranks,
with the exception of the two cases of rank 16 lattices when the dihedral group
has order 8. The latter two lattices are distinguished by $ann_M(N)$, which can
be 0 or $DD_4$. By Lemma \ref{annm(n)}, $ann_N(M) =\{ (\a,-\a)\mid \a\in
E\text{ and } h\a =-\a\}$. Therefore, $ann_M(N)\cong DD_4$ when $h$ has form
$h_{A_3}^2\oplus h_{A_1}^2$ (Theorem \ref{order4} (1) )  and  $ann_M(N)=0$ when
$h$ satisfies $h^2=-1$. Our set of rootless classes in $O(E_8)$ therefore gives
11 distinct SDC lattices, which must be the 11 types listed in \cite{GLEE8} and
which appear in column 4 of Table 2.

The main theorems \refpp{mainth1}, \refpp{mainth2}, \refpp{mainth3} of this article are now proved.  The rest of this article demonstrates new embeddings of a few of the above lattices into the Leech lattice.


\appendix

\section{Embeddings of $EE_8$ pairs in the Leech lattice}\label{embed}

As usual, $\Lambda$ denotes a copy of the Leech lattice.
\smallskip

In this appendix, we shall construct several  lattices $\mathcal{E}\cong
E_8\perp E_8$ in $\Lambda \otimes_\ZZ  \QQ$ such that $\mathcal{E}\cap \Lambda$
is an $SDC(E_8)$-lattice.   This will give relatively easy embeddings of some
rootless $EE_8$ pairs into the Leech lattice.   An account of embeddings for
all cases of $EE_8$-pairs was given in \cite{GLEE8}.

\medskip
\subsection{Order 2}

Let $\Omega$ be a  24-set and let $\mathcal{G}$ be the extended Golay code of
length $24$ indexed by $\Omega$.

For explicit calculations,  we shall use some $4\times 6$ arrays to denote the
codewords of the Golay code and the vectors in the Leech lattice. For each
codeword in $\mathcal{G}$, $0$ and $1$ are indicated by an empty and filled
space, respectively, at the corresponding positions in the
array.

\medskip

The following is a standard construction of the Leech lattice.
\begin{de}[\cite{CS,gr12}]\labtt{gen}
Let $e_{i}: =\frac{1}{\sqrt{8}}\left( 0,\dots ,4,\dots ,0\right) $ for $i\in
\Omega$. Then $(e_i,e_j)=2\delta_{i,j}$. Denote $e_{X}:=\sum_{i\in X}e_{i}$ for
$X\in \mathcal{G}$. The \textsl{standard Leech lattice} $\Lambda $ is a lattice
of rank 24 generated by the vectors:
\begin{eqnarray*}
&&\frac{1}{2}e_{X}\, ,\quad \text{where } X\text{ runs over all codewords of the Golay code }%
\mathcal{G}; \\
&&\frac{1}{4}e_{\Omega }-e_{1}\,\,; \\
&&e_{i}\pm e_{j}\,,\text{ }i,j \in \Omega.
\end{eqnarray*}
\end{de}

Let $\mathcal{D}$ be the subcode of $\mathcal{G}$ generated by
\[
\begin{split}
\mathcal{O}_1=
\begin{array}{|cc|cc|cc|}
\hline  \ *\,  &\ \, & \ *\,      &\ \   &\ \   &\ \  \\
\ *\,  & \ & \ *\,    &\ &\ &\ \\
\ *\,   & \ &\  *\,    &\ &\ &\ \\
\ *\,   & \ &\ *\,    &\ &\ &\ \\ \hline
\end{array}\ , \qquad
&
\mathcal{O}_2=
\begin{array}{|cc|cc|cc|}
\hline & \ *\,  &\ \, & \ *\,      &\ \   &\ \     \\
\ & \ *\,  & \ & \ *\,    &\ &\  \\
\ & \ *\,   & \ &\  *\,    &\ &\  \\
\ & \ *\,   & \ &\ *\,    &\ &\ \\ \hline
\end{array}\ ,\\
\mathcal{O}_3=
\begin{array}{|cc|cc|cc|}
\hline  \, *  &\ * & \ *      &\ *   &\ \   &\ \  \\
\, *  &\ * & \ *      &\ *   &\ \   &\ \  \\
\ \   & \ &\  \    &\ &\ &\ \\
\ \   & \ &\ \     &\ &\ &\ \\ \hline
\end{array}\ , \qquad
&
\mathcal{O}_4=
\begin{array}{|cc|cc|cc|}
\hline  \, *  &\, * & \, *      &\, *\,   &\ \   &\ \  \\
\, \   & \ &\  \    &\ &\ &\ \\
\, *  &\, * & \, *      &\, *   &\ \   &\ \  \\
\, \   & \ &\ \     &\ &\ &\ \\ \hline
\end{array}.
\end{split}
\]
Note that $\mathcal{D}$ is  supported at $\mathcal{O}_1\cup \mathcal{O}_2$ and
is isomorphic to $d(H_8)$,  where $H_8$ is the Hamming $[8,4,4]$-code and $d:
\ZZ_2^8 \to \ZZ_2^{16}$ is defined by $d(\a)=(\a, \a)$.

\begin{rem}
Recall that $A_1^*=\frac{1}2 A_1$, $ A_1^*/A_1\cong \ZZ_2$ and the root lattice
$E_8$ can be constructed by $A_1^8$ and $H_8$ as follows \cite{CS,gal}.

Let $\rho: (A_1^*)^8 \to (A_1^*/A_1)^8\cong \ZZ_2^8$ be the natural map. Then
$\rho^{-1} (0)=\ker \rho =A_1^8$ and $\rho^{-1} (H_8) \cong E_8$.
\end{rem}

Let
\[
A=\frac{4}{\sqrt{8}} \left\{\left.
\begin{array}{|cc|cc|cc|}
\hline  a  & b &  a & b    &\ \   &\ \  \\
c   & d & c  & d &\ &\ \\
e  & f  & e & f   &\ \   &\ \  \\
g & h &  g  & h  &\ &\ \\ \hline
\end{array} \    \right| \   a,b,c,d,e,f,g,h\in \ZZ\right\}.
\]
and denote $M=span A \cup \{ \frac{1}2 e_X\mid X\in \mathcal{D}\}$.  Then
$A\cong AA_1^8$ and $M\cong EE_8$.  Note that both $A$ and $M$ are sublattices
of $\Lambda$.

\medskip

Let
\[
\mathcal{O}=
\begin{array}{|cc|cc|cc|}
\hline  \ * & \, *\,  &\ \,      &\ \   &\ \   &\ \  \\
\ *   & \, *\,  & \  &\ &\ &\ \\
\ *  &\,  *\,  & \  &\ &\ &\ \\
\ *  &\, *\,  & \  &\ &\ &\ \\ \hline
\end{array}\ , \qquad
\hat{\mathcal{O}}=
\begin{array}{|cc|cc|cc|}
\hline \  &\ \   & \ * & \, *\,  &\ \,      &\ \     \\
&\ \   &\ *   & \, *\,  & \  &\  \\
&\ \   & \ *  &\,  *\,  & \  &\  \\
&\ \   & \ *  &\, *\,  & \  &\  \\ \hline
\end{array}
\]
and denote by $P_{\mathcal{O}}$ and $P_{\hat{\mathcal{O}}}$ the natural
projections to $\mathcal{O}$ and $\hat{\mathcal{O}}$, respectively.

\medskip

Let $E^1=P_{\mathcal{O}}(M)$ and $E^2= P_{\hat{\mathcal{O}}}(M)$. Then
$E^1\cong E^2\cong E_8$ and $E^1\perp E^2$. Moreover, $E^1\perp E^2<
\frac{1}2 \Lambda$. By identifying $E^1$ with $E^2$, we have
\[
M=\{(\a, \a)\mid \a \in E^1\}.
\]

\medskip

\noindent \textbf{Case 1:} Now let $h_1=\varepsilon_{\hat{\mathcal{O}}}$, i.e.,
$h_1$ acts as $-1$ on the basis vectors indexed by $\hat{\mathcal{O}}$ and as
$1$ on the basis vectors indexed by $\Omega\setminus \hat{\mathcal{O}}$.

Then $h_1$ acts as $-1$ on $E^2$ and fixes $E^1$ pointwise. Then
$N=h_1(M)=\{ (\a, h_1\a)\mid \a\in E^1\} <  E^1\perp E^2$ is also a diagonal
copy. In this case, $M\perp N$ and $M+ N \cong EE_8\perp EE_8$.

\medskip

\noindent \textbf{Case 2:}  Let
\[
\mathcal{O}'=  \begin{array}{|cc|cc|cc|}
\hline \ \,  &\ \, & \ *  &\ *   &\ *   &\ *   \\
 &\ & \ *  &\ *  &\ *   &\ *   \\
 &\ & &\ &\ &\ \\
 &\ & &\ &\ &\ \\ \hline
\end{array}
\]
and define  $h_2=\varepsilon_{\mathcal{O}'}$. Then $| \mathcal{O}\cap
\mathcal{O}' |=0$ and $| \hat{\mathcal{O}}\cap \mathcal{O}' |=4$. Thus, $h_2$
may be identified with $h_{A_1}^4\oplus id_{A_1}^4$ on $E^2$

and fixes $E^1$
pointwise. Let $N=h_2(M)$. Then $N\cap M \cong DD_4$ and $M+N\cong
DIH_{4}(12).$

\subsection{Order 3}
First, we recall the ternary construction of the Leech lattice $\Lambda$
\cite{CS}.  Let $\Delta$ be a 12-set and let $\TG$ be a ternary Golay code with
index set $\Delta$.

We also use the standard model for $A_2$, i.e.,
\[
A_2=\{ (a,b,c)\in \ZZ^3\mid a+b+c=0\}.
\]
Let $\gamma_0:=0$, $\gamma_1:=\frac{1}3 (1,1,-2)$ and
$\gamma_2:=\frac{1}3 (-1,-1,2)$ be elements in $A_2^*$.

Let $\mathcal{A}^i, i\in \Delta,$ be isometric copies of $A_2$ and
$\mathcal{X}:= \oplus_{i\in \Delta} \mathcal{A}_i$ an orthogonal sum of $12$
copies of $A_2$. Then the dual lattice $\mathcal{X}^* =\oplus_{i\in \Delta}
\mathcal{A}_i^*$ and $\dg {\mathcal X}$ has a natural identification with
$\FF_3^{12}$.

For each codeword $x=(x_1,\dots, x_{12})\in \TG$, let $ \gamma_x =
(\gamma_{x_1}, \dots, \gamma_{x_{12}})\in \mathcal{X}^*$ be some vector which
modulo $\mathcal X$ gives the codeword $x$. Then
\[
\mathcal{N} : = span\,  \mathcal{X} \cup \{ \gamma_x\mid x\in \TG\}
\]
is isometric to the Niemeier lattice of type $A_2^{12}$.

Let $\delta:=\frac{1}3(1,0,-1)$ be in the standard model of $A_2$
and $\hat{\delta}: =(\delta, \dots \delta)$. Then
\[
\mathcal{N}^0=\{ \a\in N\mid (\a, \hat{\delta})\in \ZZ\}
\]
is a sublattice of index 3 and has no roots.

Let $\b=(-1,1,0)\in A_2$. Then $(\b, 0,0\dots,0)+\hat{\delta}$ has  norm $4$
and the lattice $ \mathcal{N}^0 + \ZZ((\b,
0,0\dots,0)+\hat{\delta})$ is even unimodular and has no root. Hence, it is
isometric to the Leech lattice $\Lambda$ \cite[Chapter 24]{CS}.

\medskip

Next, we construct some $EE_8$ sublattices of $\mathcal{N}^0 <  \Lambda$. We
shall arrange the 12-set $\Delta$ into a $3\times 4$ array. For each codeword
in $\TG$, $0$, $1$ and $2$ are marked by a blank space and $+$ and $-$ signs,
respectively, at the corresponding positions in the array.

Let $TD$ be the subcode of $\TG$ generated by
\[
X= \begin{array}{|c|c|c|c|}
  \hline
  \ \ & + & - &\ \  \\
  \ & + & - &  \\
  \ & + & - &  \\
\hline
 \end{array},
\qquad
Y= \begin{array}{|c|c|c|c|}
  \hline
  +& + &  &  \ \ \\
  - & - & - &  \\
   &  &  +& \\
\hline
 \end{array}.
\]

Let
\[
\Omega_1= \begin{array}{|c|c|c|c|}
  \hline
  * & * &\ \,  &\ \,   \\
  \ & *&  &   \\
  \ & *  & \  &    \\
\hline
 \end{array}
\ ,\qquad
\Omega_2= \begin{array}{|c|c|c|c|}
  \hline
   &  \ \,  & *  &\ \,   \\
  * & \ & * &   \\
  \ &  & *  &    \\
\hline
 \end{array}
\]
be subsets of $\Delta$ and let $P_{\Omega_1}$ and $P_{\Omega_2}$ be the natural
projections, from $\FF_3^{\Delta}$ to $\FF_3^{\Omega_1}$, $\FF_3^{\Omega_2}$, respectively.

Then $P_{\Omega_1}(TD)$ and $P_{\Omega_2}(TD)$ are both isomorphic
to the tetracode $\mathcal{C}_4$ since they are self-orthogonal and have
dimension $2$ and length $4$.

Define a permutation $\varphi$ of $\Delta$ by
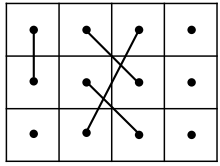
\[
\scalebox{0.7} 
{
\begin{pspicture}(0,-1.5)(4.0,1.5)
\definecolor{color0c}{rgb}{0.5019607843137255,0.5019607843137255,0.5019607843137255}
\rput(0.0,-1.5){\psgrid[gridwidth=0.028222222,subgridwidth=0.014111111,gridlabels=0.0pt,
subgriddiv=0,subgridcolor=color0c](0,0)(0,0)(4,3)}
\psdots[dotsize=0.16](0.52,1.0) \psdots[dotsize=0.16](1.52,1.0)
\psdots[dotsize=0.16](2.52,1.0) \psdots[dotsize=0.16](3.52,1.0)
\psdots[dotsize=0.16](0.52,0.02) \psdots[dotsize=0.16](1.52,0.0)
\psdots[dotsize=0.16](2.52,0.0) \psdots[dotsize=0.16](3.52,0.0)
\psdots[dotsize=0.16](0.52,-0.98) \psdots[dotsize=0.16](1.52,-0.96)
\psdots[dotsize=0.16](2.52,-1.0) \psdots[dotsize=0.16](3.52,-1.0)
\psline[linewidth=0.04cm](0.52,1.02)(0.52,0.09)
\psline[linewidth=0.04cm](1.5,1.0)(2.5,0.0)
\psline[linewidth=0.04cm](1.51,0.02)(2.5,-0.96)
\psline[linewidth=0.04cm](1.52,-0.94)(2.51,0.98)
\end{pspicture}
}
\]
Then $\varphi(\Omega_1)=\Omega_2$ and $\varphi$ induces an isomorphism between
$P_{\Omega_1}(TD)$ and $P_{\Omega_2}(TD)$.

Let
\[
B= \left\{  \left . \begin{array}{|c|c|c|c|}
  \hline
  a &\  b\, & - d &\ \  \   \\
  -a&\ c\,& -b &   \\
  \ &\  d\, & -c &    \\
\hline
 \end{array}\  \right| \  a,b,c,d\in A_2\right \} <  \mathcal{ X}
\]
and $M= span\  B \cup \{ \gamma_x\mid x\in TD\} <  \mathcal{X}^*$. Then
$B\cong AA_2^4$.

For any subset $S\subset \Delta$,  let $\tilde{P}_{S}: \mathcal{ X}^* \to
\oplus_{i\in S}\mathcal{A}_i^*$ be the natural projection.  Then
$\tilde{P}_{\Omega_1}(B)\cong \tilde{P}_{\Omega_2}(B)\cong A_2^4$.
Moreover, we have $\tilde{P}_{\Omega_1}(M)\cong
\tilde{P}_{\Omega_2}(M)\cong E_8$ since ${P}_{\Omega_1}(TD)\cong
{P}_{\Omega_2}(TD)\cong \mathcal{C}_4$, the tetracode.

Let $E^1:=\tilde{P}_{\Omega_1}(M)$ and $E^2:= \tilde{P}_{\Omega_2}(M)$.
Then $(E^1, E^2)=0$ and $E^1\perp E^2 <  \frac{1}3 \Lambda$. Note that the
permutation $\varphi$ also induces a map on $\mathcal{X}^*$ by permutating
the $\mathcal{A}_i^*$'s.  Then we have $\varphi(E^1)=E^2$ and $M=\{ (\a,
-\varphi \a)\mid \a \in E^1\} <  E^1\perp E^2$.  By identifying $E^1$ with $E^2$
using $\varphi$, we have $M=\{ (\a, -\a)\mid \a \in E^1\}\cong EE_8$.

Let $h:= h_X:= h_{A_2}^{x_1} \oplus \cdots \oplus h_{A_2}^{x_{12}}$. Note that
$h$ defines an isometry of $\mathcal{N}$ and $\Lambda$ \cite{CS,gr12}.
Moreover, $h$ acts on $E^1\perp E^2$ as $g\oplus g^{-1}$, where $g =
h_{A_2}^3\oplus id_{A_2}\in O(E_8)$. Then
\[
N= h(M) = \{ (g\a, g^{-1}\a)\mid \a \in E^1\} = \{ (\a, g\a) \mid \a\in E^1\}.
\]
In this case, $M\cap N\cong AA_2$ and $M+N\cong DIH_{6}(14)$.

\subsection{Order 5}

First  we recall a construction of the Leech lattice from $A_4^6$ \cite{CS}.

Let $S_i$, $i=1,\dots,6$, be isometric copies of $A_4$ and $S= \oplus_{i=1}^6
S_i$ an orthogonal sum of six  copies of $A_4$'s. Then the dual lattice
$S^*=\oplus_{i=1}^6 S_i^*$.

Let $\mathcal{C}$ be the subcode of $\ZZ_5^6$ generated by
\[
(1,0,1,4,4,1),\ (1,1, 0,1,4,4),\  (1,4,1,0,1,4).
\]
Then $\mathcal{C}$ is a self-dual code over $\ZZ_5$ and is a glue code
associated to the construction of $N(A_4^6)$ from $A_4^6$ \cite[Chapter
16]{CS}.

Let $a[1]:=\frac{1}5(1,1,1,1,-4), a[2]:= \frac{1}5(2,2,2,-3,-3), a[3]:=-a[2]$,
$a[4]:=-a[1]$ in $A_4^*$ and $a[0]:=0$. For each $\a=(\a_1, \dots, \a_6)\in
\mathcal{C}$, let
\[
\gamma_\a: =(a[\a_1], a[\a_2], \dots, a[\a_6]).
\]

Define
\[
\mathcal{N}: = span_{\ZZ} ( \, S \cup \{ \gamma_{\a}\mid \a \in \mathcal{C}\})  <   S^*
\]
Then $\mathcal{N}$ is isometric to the Niemeier lattice of type $A_4^6$.

Let $\eta:=\frac{1}5(2,1,0,-1, -2)$ and $\hat{\eta}:=(\eta, \eta, \eta, \eta,
\eta, \eta)$.
Then
\[
\mathcal{N}^0=\{ \a\in \mathcal{N}\mid (\a, \hat{\eta})\in \ZZ \}
\]
is an index 5 sublattice of $\mathcal{N}$ and has no roots.

Let
\[
\Lambda := span_{\ZZ} (\mathcal{N}^0\cup \{ (\b, 0,0,0,0,0) +\hat{\eta}\}),
\]
where $\beta := (-1,1,0,0,0)\in A_4$.

Then $\Lambda$ is even unimodular and has no roots. That means $\Lambda$ is
isometric to the Leech lattice \cite[Chapter 24]{CS}.

Next we shall construct some $EE_8$'s in $\Lambda$.   Let
\[
K:=\{(0,a,0,-a,-b, b)\mid a,b \in A_4\} <  S
\]
and
\[
M:=span_{\ZZ} ( \ K\cup\{ (0,a[1],0,-a[1],-a[2],a[2])\}).
\]
Then $K\cong AA_4\perp AA_4$ and $M\cong EE_8$.

Note that $$(0,1,0,-1,-2,2)= (1,0,1,4,4,1)-(1,4,1,0,1,4)\in \mathcal{C}$$ and
hence $M <  \mathcal{N}^0 < \Lambda$.

Let $P_1: S^* \to S_2^*\oplus S_6^*$  and $P_2: S^* \to S_4^*\oplus S_5^*$ be
the natural projections.

Let $E^1:=P_1(M)$ and $E^2:= P_2(M)$. Then $E^1\cong E^2\cong E_8$ and
$(E^1, E^2)=0$. By identifying  $S_2$ with $S_4$ and $S_6$ with $S_5$, we may
identify $E^1$ with $E^2$. Then, we have $M=\{(\a, -\a)\mid \a \in E^1\}$.

\medskip

Let $h:=(1, h_{A_4}, 1,  h_{A_4}^{-1}, h_{A_4}^{-2}, h_{A_4}^2)\in O((A_4^*)^6)$.
Since $(0,1,0,-1,-2,2)\in \mathcal{C}$, one can verify that
$h(\Lambda)=\Lambda$ (see \cite{Atlas} or \cite{CS}). Note that $h$ acts as
$h_{A_4}\oplus h_{A_4}^2$ on $E^1$ and as $h_{A_4}^{-1}\oplus h_{A_4}^{-2}$
on $E^2$.

Let $N:=h(M)$ and let $g:=h|_{E^1}$. Then by the identification of $E^2$ to
$E^1$,  we may identify $h|_{E^2}$ with $ g^{-1}$. Hence, we have
\[
N=h(M) =\{ (g\a ,
-g^{-1}\a)\mid \a\in E^1\} = \{ (\a,- g^{-2}\a )\mid \a \in E^1\}.
\]

In this case, $M\cap N=0$ and $M+N$ is  an SDC lattice and is isometric to
$\dih{10}{16}$.

\bigskip

{\bf Acknowledgements.} The authors thank Kavili Institute of Theoretical
Physics in Beijing for hospitality during July and August, 2010.

The first author thanks the US National Science Foundation  (DMS-0600854)
and the US National Security Agency for financial support (H98230-10-1-0201) and the Academia Sinica for hospitality during his visit in
Taipei, August 2010.

The second author thanks  National Science Council (NSC 97-2115-M-006-015-MY3)
and National Center for Theoretical Sciences, Taiwan for financial support and
University of Michigan for hospitality during a visit, March 2010.

\end{document}